\documentclass[11pt]{amsart}
\usepackage[latin1]{inputenc}
\usepackage[T1]{fontenc}
\usepackage{hyperref}
\usepackage{amssymb,calrsfs,ae}

\newtheorem*{bthm}{Main Theorem}
\newtheorem{theo}{Theorem}
\newtheorem{lemm}[theo]{Lemma}
\newtheorem{coro}[theo]{Corollary}
\newtheorem{prop}[theo]{Proposition}
\newtheorem{defi}[theo]{Definition}
\theoremstyle{definition}
\newtheorem{rema}[theo]{Remark}

\DeclareMathOperator{\Sym}{Sym}
\DeclareMathOperator{\Diff}{Diff}

\newcommand{\cD}{\mathcal{D}}

\newcommand{\cF}{\mathcal{F}}
\newcommand{\cO}{\mathcal{O}}
\newcommand{\calH}{\mathcal{H}}
\newcommand{\cL}{\mathcal{L}}
\newcommand{\setR}{\mathbb{R}}
\newcommand{\setC}{\mathbb{C}}
\newcommand{\bH}{\mathbb{H}}
\newcommand{\bI}{\mathbb{I}}
\newcommand{\bK}{\mathbb{K}}
\newcommand{\bO}{\mathbb{O}}
\renewcommand{\Im}{\mathrm{Im}}
\renewcommand{\Re}{\mathrm{Re}}
\newcommand{\Ric}{\mathrm{Ric}}
\newcommand{\Tr}{\mathrm{Tr}}

\newcommand{\ph}{{\mathrm{ph}}}

\newcommand{\CC}{\mathbb C}
\newcommand{\HH}{\mathbb H}
\newcommand{\KK}{\mathbb K}
\newcommand{\NN}{\mathbb N}
\newcommand{\OO}{\mathbb O}

\newcommand{\RR}{\mathbb R}

\newcommand{\del}{\partial}
\newcommand{\e}{\varepsilon}

\newcommand{\calA}{{\mathcal A}}
\newcommand{\calC}{{\mathcal C}}
\newcommand{\calD}{{\mathcal D}}

\newcommand{\calM}{{\mathcal M}}
\newcommand{\calO}{{\mathcal O}}

\newcommand{\calS}{{\mathcal S}}

\newcommand{\calU}{{\mathcal U}}
\newcommand{\calV}{{\mathcal V}}

\newcommand{\tr}{{\mathrm{tr}}\,}

\newcommand{\tg}{\tilde{g}}

\newcommand{\frakc}{{\mathfrak c}}

\newcommand{\KT}{{\bK\Theta}}
\newcommand{\ff}{{\mathrm{ff}}}
\newcommand{\phg}{{\mathrm{phg}}}

\title[A nonlinear Poisson transform]{A nonlinear Poisson transform \\ 
for \\ Einstein metrics on product spaces}
\author{Olivier Biquard and Rafe Mazzeo}
\address{Université de Strasbourg; Stanford University}
\email{biquard\@@math.u-strasbg.fr, mazzeo\@@math.stanford.edu}
\date{}
\thanks{O.B. supported by CNRS-USA program \# 3332. R.M. supported 
by the NSF under Grant DMS-0505709}

\begin{document}

\begin{abstract}
  We consider the Einstein deformations of the reducible rank two
  symmetric spaces of noncompact type.  If $M$ is the product of any
  two real, complex, quaternionic or octonionic hyperbolic spaces, we
  prove that the family of nearby Einstein metrics is parametrized by
  certain new geometric structures on the Furstenberg boundary of $M$.
\end{abstract}

\maketitle

\section{Introduction}
This paper is the first in a series to investigate the deformation
theory of Einstein metrics asymptotically modelled by Riemannian
globally symmetric spaces of noncompact type and of arbitrary rank. In
the special case of real hyperbolic space, and the slightly more
general setting of conformally compact asymptotically hyperbolic
manifolds, this has been the focus of extensive study over the last
fifteen years; this attention is due both to the many deep connections
with conformal geometry, cf.\ \cite{Gra00}, as well as the central
role that these `Poincaré-Einstein' spaces play in the AdS/CFT
correspondence in string theory, for which the proceedings
\cite{Biq05b} provides a good introduction. More recently, some of
this analysis has been extended to the other rank one noncompact
symmetric spaces \cite{Biq00}. Some recent advances in linear analysis
on symmetric spaces has now made it reasonable to attack this problem
in greater generality.

The rank one globally symmetric spaces of noncompact type are the
real, complex, and quaternionic hyperbolic spaces, denoted $\RR H^m$,
$\CC H^m$, $\HH H^m$, respectively, and the octonionic hyperbolic
plane $\OO H^2$.  Each has curvature bounded between two negative
constants and is diffeomorphic to an open ball $B^{n+1}$.  There is a
family of `asymptotically $\KK$ hyperbolic' ($\KK = \RR, \CC, \HH$ or
$\OO$) -- or briefly, A$\KK$H -- metrics, with asymptotics modelled on
$\KK H^m$, each of which induces a geometric structure on the sphere
at infinity, $S^n$. In the real case, this geometric structure is a
conformal class; in the complex and quaternionic case, it is a CR or
quaternionic contact structure, as described in \cite{Biq00}. There is
an octonionic contact structure on $S^{15}$, but it is rigid.
We call these boundary structures either the conformal infinity data
of the A$\KK$H metrics, or alternately, $G$-conformal structures,
where $G$ is the semisimple Lie group associated to $\KK H^m$. One
main result of \cite{Biq00} is that for each $G$ conformal structure
near the standard one on $S^n$ there is a unique A$\KK$H {\it
  Einstein} metric near to $\KK H^m$, and with that conformal infinity
data. This correspondence between A$\KK$H Einstein metrics and
$G$-conformal structures can be regarded as an asymptotic boundary
problem, and the assignment of the interior Einstein metric to the
conformal infinity data a sort of nonlinear Poisson transform. There
are other very interesting, and more subtle, problems of this type: we
mention in particular LeBrun's positive frequency conjecture,
concerning self-dual and anti-self-dual Einstein metrics in four
dimensions, which was solved by the first author in \cite{Biq02}.

It is our goal in this paper to establish a similar local deformation
theory in the first higher rank case, namely for products of the
various hyperbolic spaces listed above. There are a number of new and
interesting features not encountered in the rank one case, and the
details of the geometry and analysis are already sufficiently
complicated that it has seemed reasonable to keep this as a separate
paper. In later papers we shall treat the cases corresponding to more
general noncompact higher rank symmetric spaces. This mirrors the
recent developments for the linear analysis (for the scalar Laplacian)
\cite{MazVas02,MazVas04,MazVas05,MazVas07}.

At the roughest level, the proof proceeds exactly as in the A$\KK$H
setting, by constructing a family of approximate Einstein metrics,
parametrized by a family of boundary structures generalizing the
$G$-conformal structures, and then applying the inverse function
theorem to an appropriately gauged version of the Einstein operator.
The solution of the resulting nonlinear elliptic equation yields the
`near product hyperbolic' Einstein metric with the prescribed
conformal infinity data. There are three main new issues in carrying
this out for higher rank symmetric spaces. The first, purely geometric
in nature, involves defining the appropriate analogue of $G$-conformal
structures. This relies in turn on a choice of compactification for
each of these product hyperbolic spaces as a manifold with corners of
codimension two; the new boundary structures are defined on the
corner. However, not every one of these new boundary structures can be
extended to an asymptotically product hyperbolic metric which is also
asymptotically Einstein everywhere near infinity.  Extra hypotheses on
the boundary structure must be imposed, and even then one must solve
an extra Einstein-like equation to be able to extend this structure
from the corner to the codimension one boundary faces. The final issue
is to attain some understanding of the mapping properties of the
linearized gauged Einstein operator on these asymptotically product
hyperbolic metrics so that we can apply an inverse function theorem
argument.

In slightly more detail, if $M_j$ is a rank one hyperbolic space, $j =
1,2$, then the correct notion of conformal infinity data on $M_1 ×
M_2$ is what we call a $(G_1 × G_2)$-conformal structure on the
Furstenberg boundary $\del M_1 × \del M_2$. This is defined carefully
in §~\ref{sec:pbs}. To construct an approximate Einstein metric
associated to one of these, the first step is to extend this structure
from the corner to the full boundary $(M_1 × \del M_2) \cup (\del M_1 ×
M_2)$; this involves solving an auxiliary PDE on these hypersurface
boundaries, which is a coupled version of the Einstein equation on
each factor $M_j$. Once this has been done, we can extend the
conformal infinity data to an approximate solution of the problem,
i.e.\ a metric which induces this designated structure on the
Furstenberg boundary, and which is asymptotically Einstein uniformly
near infinity. The next step is to determine the mapping properties of
the linearized gauged Einstein operator on weighted Hölder spaces in
order to perturb this asymptotically Einstein metric to an exact one.
This is done using the geometric parametrix approach developed by the
second author and Vasy, as referenced above, which involves an
adaptation of the techniques of $N$-body scattering theory to study
elliptic theory on symmetric spaces of rank greater than one. A subtle
but important complication is that if either of the factors is
quaternionic, then the approximate Einstein metrics corresponding to
different $(G_1× G_2)$-conformal structures are not mutually
quasi-isometric, even up to diffeomorphism. Because of this, even the
function spaces vary in a nontrivial way as we vary the boundary
$(G_1× G_2)$-conformal structure. This means that we must analyze the
Green function for the linearized gauged Einstein operator not just at
the exact product metric, but for all nearby approximate Einstein
metrics. This necessitates that we carry out a parametrix construction
at `near product hyperbolic metrics', which is more complicated than
doing it just at the product hyperbolic space only.  (This difficulty
is already present in the quaternionic hyperbolic case \cite{Biq00};
there is a way to circumvent it then which unfortunately does not
generalize to this product setting, so the parametrix construction
seems unavoidable here.)

General information about the geometric analysis behind the
deformation theory for Einstein metrics, particularly in the compact
setting, can be found in \cite{Bes87}. We follow a slightly different
route developed in \cite{Biq00}. Suppose that $g$ is Einstein, i.e.\
$\Ric(g) - \lambda g = 0$ for some real number $\lambda$. If $h$ is a
sufficiently small symmetric $2$-tensor, then to make the equation $h
\mapsto \Ric(g+h) - \lambda(g+h)$ elliptic we supplement it with the
so-called Bianchi gauge condition $B^g(h) := \delta^g h - \frac12 d \tr^g
h = 0$. Equivalently, we look for solutions of the nonlinear elliptic
system
\[
N^g(h) := \Ric(g+h) - \lambda(g+h) + (\delta^{g+h})^* B^g(h) = 0.
\]
It is not hard to show that if $\lambda < 0$ and $B^g(h) \to 0$ at infinity,
then solutions of this equation correspond to Einstein metrics $g+h$
in Bianchi gauge with respect to $g$. One advantage of this gauge is
that the linearization of $N^g$ at $h=0$ is the particularly simple
operator
\[
L^g = \frac12\left(\nabla^* \nabla - 2 \overset{\circ}{R}\right),
\]
where the final term on the right is the usual action of the full
curvature tensor for $g$ on symmetric $2$-tensors. Throughout this
paper, this operator will be called the linearized gauged Einstein
operator.

Our principal result is the
\begin{bthm} Let $(M_j,g_j)$ be an A$\KK_j$H Einstein space, $j=1,2$,
  and let $\frakc$ be the product $(G_1× G_2)$-conformal structure on
  $\del M_1 × \del M_2$. Assume that $0$ is not an $L^2$ eigenvalue
  for the linearized gauged Einstein operator on either $(M_j,g_j)$,
  $j=1,2$, or $(M_1 × M_2, g_1 + g_2)$. Let $\frakc'$ be any other
  smooth $(G_1× G_2)$-conformal structure which is sufficiently
  close to $\frakc$ in the $\calC^{2,\alpha}$ norm. Assume alsgo that
  $\frakc'$ satisfies the global integrability hypothesis in Definition
  (\ref{defi:global-int}). Then there is a near product hyperbolic
  Einstein metric $g'$ with conformal infinity data $\frakc'$, and
  moreover $g'$ is unique amongst such metrics in a neighbourhood of
  the product metric $g$.
\end{bthm}

The hypothesis on the linearized gauged Einstein operator is satisfied in many 
situations, in particular for convex cocompact quotients of $\KK$ hyperbolic spaces 
and negatively curved A$\KK$H Einstein spaces.

We note also that the techniques and results of this paper apply somewhat
more generally than when $M$ is globally a product. It would not be too
difficult to define a class of manifolds and metrics with appropriate local 
product conditions near the corners and boundary faces to which this
deformation theory also applies. As a very simple example, $M$ might
be obtained by a compact topological perturbation from the product 
$M_1 × M_2$. However, we do not currently know any manifolds of this type 
which are Einstein, and so have not formulated our main result in this
greater generality.

One final comment about notation. We shall be using various classes of Hölder spaces, 
often weighted by powers of boundary defining functions. We typically write $\calC^{k,\alpha}$ 
for Hölder spaces on compact manifolds (or in any compact set), and $\Lambda^{k,\alpha}_{g}$ 
when these spaces are defined relative to some complete metric $g$ on an manifold. (In fact, 
the subscript $g$ is replaced by some moniker for a general class of complete metrics with a 
type of prescribed asymptotic geometry.)

The plan of this paper is as follows. In §~\ref{sec:r1gr} we review
the geometry of A$\KK$H hyperbolic spaces.  §~\ref{sec:rola} contains
a lengthy review of the geometric parametrix theory used to study
elliptic theory on these spaces, which is called the
$\KT$-pseudodifferential calculus; we also establish some results here
about the resolvent family of $L^g$ for such a metric.
§~\ref{sec:einst-deform-theory} reviews the analysis needed to carry
out the deformation theory of Einstein metrics in the A$\KK$H setting.
§~\ref{sec:pbs} develops the notion of $(G_1× G_2)$-conformal
structures and some geometric properties and estimates for the
corresponding asymptotically product hyperbolic metrics. The extension
of these structures to the codimension one boundary faces is the
subject of §~\ref{sec:sef}.  The parametrix construction in the near
product hyperbolic case is the topic of §~\ref{sec:gener-lapl-near},
and finally, the brief §~\ref{sec:gedt} finishes the proof of the main
theorem.

\section{Asymptotically $\bK$ hyperbolic spaces}
\label{sec:r1gr}
This section reviews the geometry of rank one symmetric spaces of
noncompact type, and of the more general class of Riemannian manifolds
which are asymptotically modelled on these.

\subsection*{Hyperbolic spaces and their conformal infinities}
The noncompact symmetric spaces of rank one are commonly called
hyperbolic spaces and written as $\bK H^m$, where $\bK = \setR$,
$\setC$, $\bH$ (the quaternions) or $\bO$ (the octonions). Note that
$\bO H^m$ exists only when $m=1,2$, and in fact $\bO H^1=\setR H^8$,
so the only new space in this last family is the $16$-dimensional
octonionic hyperbolic plane $\bO H^2$. As a homogenous space, $\bK
H^m=G/K$, where $G$ is a real semisimple Lie group and $K$ a maximal
compact subgroup; more specifically,
\begin{alignat*}{3}
  \setR H^m&=\mathrm{SO}_{1,m}/\mathrm{SO}_m,\quad & \setC H^m&=\mathrm{SU}_{1,m}/\mathrm{U}_m, \\
  \bH H^m&=\mathrm{Sp}_{1,m}/\mathrm{Sp}_1\mathrm{Sp}_m,\quad & \bO H^2&=\mathrm{F}^{-20}_4/\mathrm{Spin}_9.
\end{alignat*}
These are the noncompact duals of the corresponding projective spaces $\bK P^m$.

Throughout this paper we write
\begin{equation}
d = \dim_{\RR} \KK,
\label{defd}
\end{equation}
so that
\[
\dim_{\RR} \bK H^m = m d := n+1, 
\]
where this last equality defines $n$ in terms of $m$ and $d$. 

The polar coordinate expression for the metric on $\RR H^m$ is
\begin{equation}
g = dr^2 + \sinh^2(r)\, \gamma,
\label{eq:metrhspc}
\end{equation}
where $\gamma$ is the standard metric on $S^n$. For the analogous
expression on the other hyperbolic spaces, denote by $\eta\in
\Omega^1(S^n)\otimes\Im(\bK)$ the connection $1$-form of the Hopf bundle
\[
\begin{array}{ccc}
S^{d-1} & \longrightarrow & S^n \\ & & \downarrow \\ & & \bK P^{m-1}
\end{array}
\]
and let $\gamma$ be the pullback of the standard metric on $\bK P^{m-1}$,
regarded as a metric on the distribution $\cD=\ker \eta$. The metric on
$\bK H^m$, normalized to have sectional curvatures in $[-4,-1]$, is
given by
\begin{equation}
g=dr^2+\sinh^2(r) \, \gamma+\sinh^2(2r)\, \eta^2.
\label{eq:metohspc}
\end{equation}
The metric $\gamma$ on the distribution $\cD$ can be obtained as the limit
as $r \to \infty$ of the family of metrics $\gamma_r=4e^{-2r}g|_{TS^n_r}$ on
$S^n$; note that this limit is finite only on $\cD$, and becomes
infinite on any complementary direction.  There is no natural
`origin', so $e^{-2r}$ and $\gamma$ are defined only up to a
multiplicative factor. Indeed, once we are in the fully geometric
setting below, it is most natural to take $r$ as the distance from a
large convex hypersurface, and then we see that $\gamma$ is only
determined up to an arbitrary smooth positive factor; thus only the
conformal class $[\gamma]$ of this metric on $\cD$ is well-defined.  We
call this asymptotic data $(\cD,[\gamma])$ the \emph{conformal infinity}
of $g$.

Associated to the distribution $\cD$ on $S^n$ is the bundle $\cD\oplus
(TS^n/ \cD$) over $S^n$. The Lie bracket on sections of $\cD$ equals
$-d\eta$, and thus induces the structure of a nilpotent Lie algebra on
each fibre of this extended bundle which is isomorphic to the
$\KK$-Heisenberg algebra $\KK-\mbox{Heis}^{m-1} \cong 
\bK^{m-1}\oplus\Im(\bK)$. The metric $\gamma$ is
\emph{compatible} with $d\eta\in \Omega^2_\cD\otimes\Im(\bK)$ in the sense that the
pair $(d\eta,\gamma)$ defines a $\bK$-structure on $\cD$, i.e.\ a collection
of $d-1$ almost complex structures which are orthogonal with respect
to $\gamma$, and which satisfy the algebraic relations of the basis
elements in $\Im(\bK)$.

\subsection*{$G$-conformal structures}
The hyperbolic metrics (\ref{eq:metrhspc}) and (\ref{eq:metohspc}) are
the models for more general asymptotically hyperbolic metrics (of type
$\RR$, $\CC$, $\HH$ or $\OO$). Before defining these, however, we
first introduce terminology for the conformal infinity structures,
which will also be used later in the product case.

\begin{defi}
  Fix the $\bK$ hyperbolic space $\bK H^m=G/K$. A \emph{$G$-conformal
    structure} on an arbitrary manifold $Y^n$, $n = md-1$, is a
  codimension $d-1$ distribution $\cD \subset TY$, with a conformal
  structure $[\gamma]$ on the fibres of $\cD$, such that the induced
  nilpotent Lie algebra $\cD\oplus TY/ \cD$ is isomorphic at each point to
  the $\bK$-Heisenberg algebra, and any metric $\gamma \in [\gamma]$ is
  compatible with the $\KK$ structure on $\cD$. We say that the distribution $\cD$ is of
  $\bK$-contact type.
\end{defi}
This definition unifies several cases:
\begin{itemize}
\item when $\bK = \RR$, $\cD$ is the entire tangent space 
  and $[\gamma]$ is a conformal structure in the usual sense;
\item when $\bK = \CC$, $\cD$ is a contact distribution in the ordinary sense; 
  if $\eta$ is a contact $1$-form which defines $\cD$, then the compatibility 
  of a metric $\gamma$ on $\cD$ with $d\eta$ means that on $\cD$ one has 
  $d\eta(\cdot,\cdot)= \gamma(J\cdot,\cdot)$ for some almost complex   structure $J$
  which is orthogonal with respect to any $\gamma \in [\gamma]$; this is simply 
  an almost CR structure on $\cD$;
\item when $\bK = \HH$, $\cD$ is a ``quaternionic contact structure'' 
  as defined and studied in \cite{Biq00}; it turns out that 
  the conformal class $[\gamma]$ is completely determined by $d\eta$;
\item finally, the octonionic case is rigid; $\cD$ is automatically 
  locally isomorphic to the standard distribution on the sphere $S^{15}$, 
  and $[\gamma]$ is determined completely by $d\eta$.
\end{itemize}

It is important to note here that in the quaternionic case, even
though each of the tangent nilpotent Lie algebras is isomorphic to the
standard quaternion Heisenberg algebra, the distribution $\cD$ is not
locally diffeomorphic to the model structure on the sphere (unless it
is standard everywhere). In other words, there is no direct analogue
of Darboux's theorem for quaternionic contact structures, and the
infinitesimal equivalence of these structures at each point does not
imply their local equivalence.

\subsection*{A$\bK$H metrics}
Let $M$ be a manifold with boundary $Y = \del M$ admitting a
$G$-contact structure (associated to $\KK H^m = G/K$). We give two
equivalent definitions of the class of complete metrics on the
interior of $M$ which induce a $G$-conformal structure on $Y$.

The first mimics the polar coordinate definition of the model case.
\begin{defi}
  A metric $g$ on $M^{n+1}$ is called \emph{asymptotically $\bK$
    hyperbolic} (or A$\bK$H for short) if the following conditions are
  satisfied: there is a neighbourhood $\calU$ of $Y$ in $M$, a
  diffeomorphism identifying $\calU$ with $(1,\infty)_r × Y$, and a
  $G$-conformal structure $(\cD,[\gamma])$ on $Y$, such that, fixing a
  representative $(\eta,\gamma)$ of the $G$-conformal structure and defining
\[
g_0(\gamma,\eta) = dr^2+\sinh^2(r)\, \gamma + \sinh^2(2r)\, \eta^2
\]
in $\calU$, we have
\[
g= g_0(\gamma,\eta) + k, 
\]
where $k$ is in the weighted geometric Hölder space $e^{-\nu r} \Lambda^{2,\alpha}$ 
for some $\nu>0$.  (The derivatives and norms are taken with respect to $g_0(\gamma,\eta)$.)
\end{defi}
This definition does not depend on the choice of $(\eta,\gamma)$ in the
conformal class since, replacing $(\eta,\gamma)$ by $(f\eta,f\gamma)$ for some $f
\in \calC^\infty(\overline{M})$, $f > 0$, changes the model (up to
diffeomorphism) by an error which is $\calO(e^{-r})$. The pair
$([\gamma],\eta)$ (or more properly, the triple $(\cD, [\gamma], \eta)$), is
called the conformal infinity of $g$.

The alternate definition simply replaces the radial variable $r$ by $x
= e^{-r}$, which is a defining function for $Y$ in $M$ (recall, this
means that $x \geq 0$ in $M$, $x = 0$ only on $Y$ and $dx \neq 0$ there);
thus
\begin{equation}
g_0(\gamma,\eta) = \frac{dx^2 + \gamma}{x^2} + \frac{\eta^2}{x^4},
\qquad g = g_0(\gamma,\eta) + k,
\label{eq:thmet}
\end{equation}
where $k\in x^\nu \Lambda^{2,\alpha}$. (Again, norms and derivatives are with
respect to $g_0(\gamma,\eta)$.) This will be more useful from our point of
view since the boundary $Y$ appears explicitly as the hypersurface
$\{x=0\}$. Near the boundary, the volume form of $g$ has the form
\begin{equation}
dV_g = x^{-n-d}\, dx \, dA_Y(x)
\label{eq:volform}
\end{equation}
where $dA_Y(x)$ is a family of volume forms on $Y$ depending smoothly
on $x$.

A straightforward calculation, cf.\ \cite{Biq00}, shows that an
A$\KK$H metric has curvature tensor which is asymptotic to that of
$\bK H^m$ to order $\calO(e^{-\nu r}) = \calO(x^\nu)$.

The complex hyperbolic metric on the unit ball $B$ in $\CC^m$ has a
slightly different form in standard Euclidean coordinates, and it is
worth explaining the difference.  This metric has Kähler form
\[
- \frac{\partial \overline{\partial} \rho}{\rho} + \frac{\partial \rho}{\rho}
\land \frac{\overline{\partial} \rho}{\rho},
\]
where $\rho = \frac12 (1-|z|^2)$ is a defining function for $\partial B$.  As
a Hermitian metric, the first term blows up only like $1/\rho$; its
leading coefficient is the Levi form, which is positive definite on
$\cD$.  The second term, which blows up at the faster rate $1/\rho^2$,
vanishes on (the radial extension of) $\cD$, and is positive on the
directions spanned by $\partial \rho$ and $\overline{\partial} \rho$, or equivalently,
on the span of $\del_\rho$ and $i \del_\rho$. There are analogous
expressions for the Bergman and Kähler-Einstein metrics on any
strictly pseudoconvex domain.

The obvious discrepancy with (\ref{eq:thmet}) is resolved by setting $x = \sqrt{\rho}$.  
This accords with the fact that the geodesic distance function $r$ for the hyperbolic metric 
is comparable to $-\frac12 \log \rho$ rather than $-\log \rho$.  More bluntly, the standard
$\calC^\infty$ structure on the closure of the Euclidean ball (or any strictly pseudoconvex domain) 
induced from its inclusion in $\CC^m$ is not quite the right one for our purposes.

\section{Linear elliptic theory on asymptotically hyperbolic spaces}
\label{sec:rola}

We now describe the structure of the Green function and mapping
properties for the linearized gauged Einstein operator on an
asymptotically $\KK$ hyperbolic space.

There are several ways to approach linear elliptic problems of this
type. Because the underlying geometric structure is asymptotically
rank one, certain features of the operators in question are dominated
by their radial behaviour, which is one-dimensional and hence more
readily tractable.  Using this, the first author \cite{Biq00} carried
out a detailed ODE analysis for the radial part of the relevant
operators on each $\KK$ hyperbolic space to capture the decay of the
corresponding Green functions, from which the required mapping
properties can be deduced. In principle, the same general ideas should
work for higher rank geometries, but the radial parts of these
operators are then multi-dimensional and must be studied differently.
This will be done later in this paper via the more general techniques
of geometric microlocal analysis. We review these methods in the setting of
asymptotically $\KK$ hyperbolic geometry, even though simpler methods
are available there, because this is a good warm-up for the construction
in the product case below, but also since we require certain more subtle
estimates on the resolvent family here which are used in the product
analysis of §7.

Before embarking on all of this, let us say a few words about the
general strategy. Local elliptic theory, or elliptic theory on compact
manifolds, can be developed entirely via Schauder estimates.  A
satisfactory understanding of the global mapping properties for a
Laplace-type operator $L$ on a complete noncompact manifold requires
not only this local theory but also some information about `far-field'
effects. Roughly speaking, one needs estimates at infinity for
solutions of $Lu = f$, even when $f \in \calC^\infty_0$. Supposing for
simplicity that $L$ is actually invertible on $L^2$, then its inverse
is represented by an integral operator $f \mapsto u(z) = \int
G(z,z')f(z')\, dz'$. This integral kernel, $G$, is a distribution on
$M × M$, and is called the Green function for $L$.  Its structure near
the diagonal $\{z=z'\}$ is exactly the same as in the compact case,
but the interesting mapping properties of $L$ are determined by its
asymptotics as $(z,z')\to \infty$ in any direction in $M × M$. For example,
for nonlinear problems one must usually understand the invertibility
of $L$ on Hölder spaces rather than Sobolev spaces. It is not so easy
to deduce Hölder boundedness from $L^2$ boundedness on a noncompact
space, but this follows from pointwise estimates for the off-diagonal
asymptotics of $G$. Thus it is a fundamental goal to determine these
asymptotics, at least for special classes of complete spaces.

The existence of $G$ may be known by abstract or indirect methods,
e.g.\ using Hilbert space theory and a Bochner-type argument, but
these usually give little information about the asymptotic structure
of this Schwartz kernel. A parametrix for $L$ is an approximation to
$G$, or slightly more generally (in case $L$ is only Fredholm) an
operator which inverts $L$ up to compact errors (though note that
compactness depends on the function spaces on which these act). A
parametrix construction produces operators for which we have good
pointwise control on the Schwartz kernels. In these special geometric
settings, such a construction proceeds by solving a sequence of model
problems, and using their solutions to construct successively better
approximations to the true (putative) Green function. To organize this
information in a useful manner, we begin by defining a certain
compactification $\widetilde{M^2}$ of $M × M$. This compactification
is a manifold with corners, and the model problems appear as the
induced operators on its different boundary hypersurfaces. The fine
pointwise structure of a parametrix is encoded in the statement that
it is a conormal (or even better, a polyhomogeneous conormal)
distribution on $\widetilde{M^2}$. A posteriori one also deduces this
same regularity structure for the Green function itself. The
technicalities of this construction involve defining a
pseudodifferential calculus on $M$ which is large enough so that one
may carry out some sort of parametrix construction for `fully
elliptic' operators. These operators are characterized by the
regularity properties of their Schwartz kernels on $\widetilde{M^2}$,
and the main work consists in verifying the usual properties, i.e.,
composition, boundedness, etc. This can be done in a number of
situations under fairly strong assumptions on the ambient geometry,
all satisfied in this asymptotically hyperbolic (or higher rank
symmetric) setting.  In this section we explain how to carry this out
for the asymptotically $\KK$ hyperbolic geometries.

We assume in this section that $L$ is a generalized Laplacian acting
between sections of bundles $E$ and $F$ over $M$ which are associated
to the bundle of orthonormal frames (or else as a parallel subbundle
when the holonomy of the metric is reduced) via a representation of
the orthogonal group; such bundles come equipped with the Levi-Civita
connection $\nabla$. For example, these could be parallel subbundles of
the full tensor bundle. A \emph{geometric differential operator} is a
linear combination of powers $\nabla^k$ of the covariant derivative on one
of these bundles, with coefficients determined by the metric and
curvature tensor. In particular a \emph{generalized Laplacian} on a
geometric bundle $E$ is an operator of the form
\[
L = \nabla^* \nabla + R,
\]
where $R$ is a symmetric endomorphism on $E$ constructed from the curvature tensor associated to $\nabla$. 
(However, all that we say here adapts easily to first order Dirac-type operators, and to many other 
operators besides.) In this section we describe the construction of a parametrix for $L$, and explain 
how it gives information about the Green function and mapping properties on various function spaces. 
At the end of this section we also collect some additional facts about the Schwartz kernel of the resolvent 
family $R(\lambda) = (L - \lambda)^{-1}$; this is an important ingredient in the construction of the
Green function in the product case. 

We conclude with some historical remarks. The analysis of elliptic
uniformly degenerate operators (which is the case $\bK = \RR$)
appeared in \cite{Maz88} and \cite{Maz91}, but see the monograph
\cite{Mazzeo-EMS} for a more comprehensive treatment and many
applications. The complex case was developed initially by Epstein,
Melrose and Mendoza \cite{EpsMelMen91}; further development and
ramifications of this theory are contained in the unpublished
manuscript \cite{EM-tubes}. The quaternionic and octonionic cases have
not been written down explicitly before, though the way to do so had
certainly been clear from \cite{EpsMelMen91}. These various
pseudodifferential calculi are quite similar to one another, with only
minor and obvious modifications needed between them. One small issue,
which presents only minor difficulties, is the lack of a Darboux
theorem when $\bK = \bH$. Using language introduced below, the point
is that these constructions depend on the infinitesimal, rather than
the local, identifications of manifolds with $\KT$ structure and the
$\bK$-Heisenberg models.

\subsection*{$\KT$ structures} 
We first describe the notion of a $\KT$ structure on a manifold with boundary $M$, 
and its ancillaries: the $\KT$ tangent, cotangent and tensor bundles, and the
classes of $\KT$ metrics and $\KT$ differential operators, the latter of which 
contains all geometric elliptic operators for any $\bK\Theta$ metric as elliptic elements.  

Suppose that $Y = \del M$ carries a distribution $\cD$ of type $\bK$;
let $([\gamma], \eta)$ be any associated conformal infinity. Choosing an
identification of a neighbourhood $\calU$ of $Y$ in $M$ as a product
$Y × [0,1)_x$, we extend this data to $\calU$, and hence write down
the model A$\bK$H metric $g_0 = g_0(\gamma,\eta)$ as in (\ref{eq:thmet}).
Now define the space $\calV_\KT$ of all smooth vector fields $V$ on
the (closed) manifold $M$ such that $g_0(V,V)$ is smooth (in $\calU$)
up to $x=0$. This is independent of $\gamma$ and $\eta$, but depends on
$\cD$ (and the $1$-jet of its extension to the interior).

It is helpful to write this out in a local frame. First choose a local
frame $Y_1, \ldots, Y_\ell$, $\ell = m(d-1)$, for $\cD$ and another set of
independent vector fields $Z_1, \ldots, Z_r$, $r = d-1$, which are
complementary to $\cD$ at each point and tangent to each $Y × \{x\}$;
the vector field $\del_x$ completes this to a full basis of sections
of $TM$. Then $V \in \calV_{\KT}$ if and only if
\begin{equation}
V = a\, x\del_x + \sum_{j=1}^\ell b_j\, xY_j + \sum_{k=1}^r c_k\, x^2 Z_k,
\label{eq:vkth}
\end{equation}
where $a, b_j, c_k$ are all $\calC^\infty$ up to the boundary. Hence $\calV_{\KT}$
is the span over $\calC^\infty(\overline{M})$ of $\{x\del_x, xY_1, \ldots, xY_\ell, x^2 Z_1, 
\ldots, x^2 Z_r\}$. 

The terminology `$\Theta$-structure' comes from \cite{EpsMelMen91}, where
$\Theta$ denotes a nonvanishing section of $T^*_{\del M}M \otimes \Im(\KK)$, the
pullback of which to $T^* \del M$ equals $\eta$. Hereafter, we let $\Theta$
denote not only this form on $\del M$, but also some choice of smooth
extension to the interior. All of the notions here can be defined in
terms of this form, or equivalently, the corresponding oriented
$(d-1)$-plane bundle $S$ of $1$-forms on $M$ at $\del M$. Thus, for
example, elements of $\calV_{\KT}$ are also characterized as vector
fields which are smooth on $\overline{M}$ and which satisfy $\Theta(V) =
\calO(x^2)$.

Note that $\calV_{\bK\Theta}$ is closed under Lie bracket. Next, there is
a vector bundle, ${}^{\KT}TM$, for which $\calV_{\KT}$ is the {\it
  entire} space of smooth sections. The fibres are defined by
\[
{}^{\KT}T_pM = \calV_{\KT}/ {\mathcal I}_p\calV_{\KT},
\]
where ${\mathcal I}_p$ denotes the space of smooth functions on $M$
vanishing at $p$. (For a more prosaic definition, we take the sections
$x\del_x, xY_i, x^2 Z_j$ as a local basis of sections of ${}^\KT TM$.)
This bundle is naturally isomorphic to $TM$ over the interior, but the
natural bundle map
\[
{}^{\KT}TM \longrightarrow TM,
\]
defined via evaluation, $V \mapsto V(p)$, is the zero map when $p \in
Y$. The subbundle over $\del M$ spanned by $\{xY_i, x^2 Z_j\}$ can
also shown to have an invariant definition, and we denote (with a
slight abuse of notation) by ${}^\KT T\del M$. Third, when $p \in Y$,
the subspace ${\mathcal I}_p \calV_\KT$ is an ideal in $\calV_\KT$
with respect to bracket of vector fields; hence for such $p$, ${}^\KT
T_pM$ is a Lie algebra, where the Lie bracket of two elements given as
the equivalence class of the vector field bracket of representatives
of the two individual classes, and as such is isomorphic to the
solvable homogeneous extension $\bK\calS^m$ of
$\bK\mathrm{Heis}^{m-1}$, the $\bK$-Heisenberg algebra. Note that
${}^\KT T_p\del M$ is a nilpotent subalgebra, isomorphic to
$\bK\mathrm{Heis}^{m-1}$ itself.

In the simplest case, when $\bK = \RR$, this space of vector fields is
usually called the space of uniformly degenerate vector fields,
denoted $\calV_0$, and consists of all smooth vector fields on
${\overline M}$ vanishing at $Y$. There is no form $\Theta$ now, so to
avoid complicating the presentation we shall mostly discuss only the
other cases, save for a few passing comments about the real case. The
other familiar case is when $M$ is a strictly pseudoconvex domain in
$\CC^m$; the $\CC\Theta$ structure is the CR structure on the boundary and
the canonical Bergman or Kähler-Einstein metrics are $\CC\Theta$ metrics
(admittedly only polyhomogeneous rather than $\calC^\infty$). As explained
earlier, one needs to take $x$ as the square root of the Euclidean
distance to the boundary in order to fit this into the present
framework.

The dual of the $\KT$ tangent bundle is denoted ${}^{\KT}T^*M$. Note
that smooth sections of this $\KT$ cotangent bundle are singular in
the ordinary sense: in terms of the dual basis of one-forms $dx$,
$Y_i^*$ and $Z_j^*$,
\[
\calC^\infty(M;{}^\KT T^*M) \ni \omega = a\, \frac{dx}{x} + \sum b_i\, \frac{Y_i^*}{x} + 
\sum c_j\, \frac{Z_j^*}{x^2},
\] 
where $a, b_i, c_j \in \calC^\infty(M)$. Similar remarks apply to all other
tensor bundles too. Note in particular that an A$\bK$H metric $g$ is a
section of $S^2({}^{\KT}T^*M)$ which is positive definite on
${}^{\KT}TM$.  If $E$ is any bundle constructed functorially from
$TM$, then applying the same functorial operations to ${}^\KT TM$
yields a bundle which we denote ${}^{\KT}E$.

\begin{defi}
  Let $M$ be a compact manifold with boundary and $\cD$ a distribution
  of type $\bK$ on $Y = \del M$.  The space $\Diff^*_{\KT}(M)$ of
  $\KT$ differential operators on $M$ consists of all operators which
  can be locally expressed as a finite sum of products of elements of
  $\calV_{\KT}$. If $E,F$ are any vector bundles over $M$, then a
  $\KT$ operator acting between sections of $E$ and $F$ is one which
  has this form with respect to any local trivialization.
\end{defi}
It is important for $\calV_{\KT}$ to be closed under Lie bracket for
this space of operators to be well-defined.

In our applications, the bundles $E$ and $F$ are geometric bundles, and the
operator is a $\KT$ differential operator between ${}^\KT E$ and ${}^\KT F$.

\begin{theo} Let $g$ be any A$\bK$H metric on $M$, and $L$ a geometric
  elliptic operator of order $\mu$ between sections of two geometric
  bundles $E$ and $F$. Then $L \in \Diff^\mu_{\KT}(M; {}^\KT E,{}^\KT
  F)$.
\end{theo}

This result is tautological once one checks that the Levi-Civita
connection $\nabla$ satisfies
\[
\nabla: \calC^\infty(M;E) \longrightarrow \calC^\infty(M; E \otimes {}^{\KT}T^*M).
\]
We leave details to the reader. 

There is a principal symbol mapping for $\calV_{\KT}$ operators,
defined formally replacing $x\del_x$, $xY_i$ and $x^2 Z_j$,
respectively, by linear coordinates $\xi, \eta_i, \zeta_j$. Thus,
\begin{equation*}
\begin{split}
\Diff^\mu_{\KT}(M;E,F) \ni P = \sum_{j+|\alpha|+|\beta| \leq \mu}
a_{j\alpha\beta}(z) (x\del_x)^j (xY)^\alpha (x^2 Z)^\beta \\
\longmapsto  \qquad {}^{\KT}\sigma_\mu(P)(z;\xi,\eta,\zeta) \quad = \quad 
\sum_{j+|\alpha|+|\beta| =\mu}a_{j\alpha\beta}(z) \xi^j \eta^\alpha \zeta^\beta.
\end{split}
\end{equation*}
The usual calculation shows that this is a well-defined smooth
function on ${}^{\KT}T^*M$ (with values in $\mathrm{Hom}(E,F)$),
homogeneous of degree $\mu$ on the fibres.

\begin{defi} The operator $P \in \Diff^\mu_{\KT}(M;E,F)$ is called
  ($\KT$) elliptic if ${}^{\KT}\sigma_\mu (P)(z;\xi,\eta,\zeta)$ is an invertible
  endomorphism whenever $(\xi,\eta,\zeta) \neq 0$.
\end{defi}

\subsection*{Parabolic dilations and model operators}
The key to the analysis of $\KT$ operators is their approximate
dilation invariance.  More precisely, for any $p \in Y$, one may define
an equivalence class of dilations based at $p$. When $\bK = \RR$,
these are ordinary radial dilations, but in the other cases the
dilations are `parabolic'. Using these we can define for any $P \in
\Diff^*_{\KT}$ the normal operator $N_p(P)$; this is a finite
dimensional reduction in that it is a left-invariant operator on the
solvable group $\bK\calS^m$, depending parametrically on $p \in Y$. Its
invertibility (for all $p$) is the other key hypothesis, besides
symbol ellipticity, needed to prove that $P$ is Fredholm.

We begin by defining these families of dilations. The situation is
simplest when $\bK = \RR$; in this case, choose a diffeomorphism of a
neighbourhood of $p$ in $M$ with a half-ball around the origin in the
half-space $\RR^n_{+} = \{(s,u): s \geq 0, u \in \RR^{n-1}\}$. Now use
this identification and the ordinary dilation operator $M_\delta: (s,u)
\mapsto (\delta s, \delta u)$ to define the sequence of pushforwards of the
vector field $V \in \calV_0$:
\begin{equation}
\lim_{\delta \to 0} (M_\delta^{-1})_* V = N_p(V).
\label{eq:nopvf}
\end{equation}
{}From the local coordinate description of $V$, we see readily that
$N_p(V)$ is a left-invariant operator defined on $\RR^n_+ \cong
\RR\calS^n$. More generally, for any uniformly degenerate differential
operator $P$, define
\begin{equation}
\lim_{\delta \to 0} (M_\delta^{-1})_* P = N_p(P);
\label{eq:nopop}
\end{equation}
This is in the universal enveloping algebra of $\RR\calS^n$, and is
well-defined up to the action of an element $A \in \RR\calS^n$.

In the other cases we begin by recalling the parabolic dilation
structure on $\KK\calS^m$.  To define this, recall that $\KK\calS =
\RR^+ \ltimes \KK\mathrm{Heis}^{m-1}$ (this is just the $A \cdot N$ part
of the $G = KAN$ decomposition). Choose a system of coordinates
$(s,\sigma,u)$ where $s > 0$, $\sigma \in \RR^{d-1}=\Im\, \KK$ and $u \in
\RR^{d(m-1)}=\KK^{m-1}$ so that (in $\KK$ coordinates)
\begin{equation}
\Theta_0 = d\sigma + \frac{1}{2} \Im(du \cdot \bar{u})\label{eq:Theta0}
\end{equation}
defines the standard $\KK$ contact structure on
$\KK\mathrm{Heis}^{m-1}$.  Thus, 
\begin{equation}
s \partial_s, s \big(\partial_{u_i} - \tfrac{1}{2} \Im(\partial_{u_i} \bar{u})\big) , 
s^2 \partial_{\sigma_j}\label{eq:q-left-inv}
\end{equation}
is a basis of left-invariant vector fields, where in $\Im(\partial_{u_i}
\bar{u})$ we identify the vector $\partial_{u_i}$ with a vector with $\KK$
coordinates in $\KK^{m-1}$, and the imaginary quaternions with
vertical vectors $\partial_{\sigma_j}$.  The dilation is then given by
\[
M_\delta(s,\sigma,u) = (\delta s, \delta^2 \sigma, \delta u).
\] 
(Note that these vector fields are homogeneous of degree $0$ with
respect to $M_\delta$.)  When $\bK = \CC$ or $\OO$, we can choose a
diffeomorphism as before which identifies a neighbourhood of $p$ in
$M$ with a neighbourhood of $0$ in $\bK\calS$, which carries the
distribution $\cD$ to the model distribution on
$\bK\mathrm{Heis}^{m-1}$; the model $\Theta_0$ above is a suitable choice
for $\Theta$ on $M$. In the complex case this uses the Darboux theorem,
while in the octonion case this follows from the local rigidity of
octonion contact structures (so that $(Y,\cD)$ is locally identified
with the model geometry).  In terms of this identification, we define
$N_p(P)$ by the same formula as above, arriving at an operator which
is left-invariant on $\KK\calS$ and well-defined up to translation by
an element of this group. In the last case, $\bK = \HH$, one has no
longer the Darboux theorem or rigidity, but the following Lemma is
proved in the Appendix.
\begin{lemm}\label{lem:coord-q}
  For any quaternionic contact structure on $Y^{4m-1}$, and any point
  $p\in Y$, there exist local coordinates $(\sigma,u)\in \Im(\bH)×\bH^{m-1}$,
  such that the quaternionic distribution is given by the kernel of a
  1-form $\Theta$ with values in $\Im(\bH)$, and the difference with the
  standard form $\Theta_0$ of the Heisenberg group near the origin
  $(\sigma=0,u=0)$ satisfies the estimate $|\Theta-\Theta_0|=O(|u|^2+|\sigma|)$.
  The two nilpotent algebra structures coincide at the point $p$.
\end{lemm}
This result constructs a diffeomorphism from a neighbourhood of $p$
to a neighbourhood of $0$ in $\HH\calS$ so that the distributions
agree at the origin. Letting $\{Y_i\}$ and $\{Y_i^0\}$ be local frames
for $\cD$ and the model distribution (in $\HH\mathrm{Heis}^{m-1}$,
extended to $\HH\calS^m$ and then transfered to this neighbourhood),
then clearly $Y_i = Y_i^0 + O(|u|^2+|\sigma|)$. It follows that the limit
of the parabolic dilations of $P$ is still a left-invariant operator
on $\HH\calS$.

To express this more concretely, fix a boundary defining function and
smooth vector fields $Y_i$, $i = 1, \ldots, d(m-1)$, $Z_j$, $j=1, \ldots,
d-1$, such that the $Y_i$ span the extension of $\cD$ and the $Z_j$
span a subspace complementary to $\cD$ at each point.  Using the
obvious multi-index notation, write
\[
P = \sum_{j + |\alpha| + |\beta| \leq m} a_{j\alpha\beta}(w)(x\del_x)^j(xY)^\alpha (x^2 Z)^\beta,
\]
where the coefficients are assumed to be $\calC^\infty$ (up to the
boundary). The values of $x\del_x$, $xY_i$ and $x^2Z_j$ at $p$ fix an
isomorphism of ${}^{\KT} T^*_pM$ and $\bK\calS$, and
\[
N_p(P) = \sum_{j+|\alpha|+|\beta|\leq m} a_{j\alpha\beta}(p)(s\del_s)^j(sY^0)^\alpha(s^2Z^0)^\beta
\]
where $s \in \RR^+$, and $\{Y_1^0, \ldots, Y_{d(m-1)}^0, Z_1^0, \ldots,
Z_{d-1}^0\}$ are a fixed basis of left-invariant vector fields on
$\bK\mathrm{Heis}^{m-1}$.

The following result is well-known in the real and complex cases. It
is obvious in the octonionic case, and a direct consequence of Lemma
\ref{lem:coord-q} in the quaternionic case.
\begin{prop}
  Let $g$ be an A$\bK$H metric on $M$, and $L$ be a
  generalized Laplace operator on $M$. Then, at each point of $\partial M$, the normal
  operator of $L$ is the corresponding operator on the hyperbolic
  space $\bK H^m$. In particular, it does not depend on the point of $\partial M$.
\label{pr:idno}
\end{prop}

For $P$ as above, there is a simpler family of model ordinary
differential operators on $\RR^+$ called the indicial family, defined
by the expression
\[
I_p(P) = \sum_{j \leq m} a_{j00}(p)(s\del_s)^j.
\]
The coefficients are endomorphisms of $E_p$. Since this is a constant
coefficient Fuchsian operator, it is equivalent by Mellin transform to
multiplication by a (matrix-valued) polynomial
\[
I_p(P;\zeta) = \sum_{j \leq m} a_{j00}(p)(\zeta)^j.
\]
A number $\zeta \in \CC$, is called an indicial root if $I_p(P;\zeta)$ is singular. This
is equivalent to the requirement that
\[
P(x^\zeta v(y)) = \calO(x^{\zeta+1}) \qquad \forall\ v \in \calC^\infty(Y).
\]
These indicial roots are fundamental invariants of $P$.

\subsection*{Blowups and the $\KT$ double space}
There is a more sophisticated way to interpret the parabolic
dilations, leading to a more obviously invariant definition of normal
operators. The idea is to introduce a resolution, or blowup, of the
product space $M × M$ which reflects the scaling invariance properties
of $\KT$ differential operators near the boundary. This provides the
means to define the $\KT$ pseudodifferential operators. We describe
this now.

As usual, we start with the simplest case $\KK = \RR$. The
distributional Schwartz kernels of pseudodifferential operators are
singular along the diagonal in $M × M$. Unfortunately, this diagonal
intersects the corner $\del M × \del M$, making it difficult to
describe the precise structure of its singularity near this
intersection. To remedy this we introduce a new space
\[
M^2_0 = [M × M; \mathrm{diag}_{\del M × \del M}],
\]
where this notation on the right indicates that we blow up $M × M$ at
the boundary of the diagonal.  This amounts to replacing this
submanifold by the space of inward pointing unit normal vectors; the
space $M^2_0$ is endowed with the smallest $\calC^\infty$ structure
containing the lifts of all smooth functions on $M× M$ and polar
coordinates around $\mathrm{diag}_{\del M}$. Thus $M^2_0$ has three
hypersurface boundaries, $B_{10}$ and $B_{01}$, the left and right
faces, which are the ones lifted from the two hypersurface boundaries
in $M× M$, and the new front face $B_{11}$ created in this blowup,
which is often also denoted $\ff$.  The blowdown map $\beta: M^2_0
\longrightarrow M^2$ is a smooth mapping of manifold with corners.

The front face $B_{11}$ fibres over $\mathrm{diag}_{\del M}$, with fibre
at $p \in Y$ the set of unit inner normal vectors at that point; this
is a quarter-sphere, the interior of which carries a natural
projective structure. Let $(x,y)$ and $(x',y')$ denote coordinates on
the two copies of $M$ in $M^2$; we are blowing up the submanifold
$x=x' = 0$, $y=y'$, and so it is legitimate to introduce the new
singular coordinate system $s = x/x'$, $u = (y-y')/x'$, $x'$, $y'$.
The coordinates $(s,u)$ are then projective coordinates on this
quarter sphere. The normal operator of $P$ is the restriction to the
fibres of $B_{11}$ of the lift of $P$ from the left factor of $M$ to
$M^2$ and then to $M^2_0$.  Thus, as in the previous definition, each
$N_p(P)$ acts on a half-space $\RR^n_+ = \RR^+ \ltimes \RR^{n-1}$.
The underlying dilation structure is implicit here since we are taking
the normal blowup, which involves ordinary homothetic scaling in the
tangent spaces.

There is a similar development for the other cases, but the normal
blowup of the boundary of the diagonal must be replaced by a blowup of
this submanifold which respects the underlying parabolic dilation
structure. Now, instead of ordinary spherical normal vectors, we use
equivalence classes of paths converging to $p$, where the equivalence
relationship is governed by the form $\KT$.

First define $I_Y$ to consist of the smooth functions on $M$ vanishing
on $Y$ and, recalling the bundle $S \subset T_{\del M}^*M$ determined by
$\Theta$, let $I_S$ denote the subset of those elements $f \in I_Y$ such
that $\left. df \right|_{Y}$ are sections of $S$. Next, fix $p \in Y$
and define the set of $S$-parabolic curves at $p$ to consist of those
smooth functions $\gamma: [0,1]\to M$, $\gamma(0) = p$ and $f(\gamma(t)) =
\calO(t^2)$ for all $f \in I_S$. We define an equivalence relation on
such curves:
\begin{equation*}
\begin{split}
  \gamma_1 \sim \gamma_2 \Leftrightarrow f(\gamma_1(t))-f(\gamma_2(t)) = \calO(t^2)
  \qquad \forall f \in I_Y; \\
  f(\gamma_1(t))-f(\gamma_2(t)) = \calO(t^3)\qquad \forall f \in I_S.
\end{split}
\end{equation*}
The space of equivalence classes is the set of inward-pointing
$S$-parabolic normal vectors to $Y$ at $p$, which we denote $\calS_p$,
and these fit together to form a bundle $\calS$ over $Y$. Each
$\calS_p$ has a natural $\RR^+$ and additive structure, defined by $\delta
\, [\gamma(t)] = [\gamma(\delta t)]$ and $[\gamma] = [\gamma_1] + [\gamma_2]$ if $f(\gamma(t))-
(f(\gamma_1(t))+f(\gamma_2(t)) = \calO(t^2)$ for all $f \in I_Y$, and
$\calO(t^3)$ for all $f \in I_S$, respectively.  These do {\it not}
define a linear structure, however, since the scalar action does not
distribute over addition.  However, directly from a local coordinate
calculation one finds that $\calS_p\cong \bK\calS$, and scalar
multiplication corresponds to parabolic dilation.

\subsection*{Conormal distributions}
We make a small diversion from the main thread of this section to
define various classes of conormal distributions on a manifold with
corners; these will appear in many places below. Let $X$ be a manifold
with corners and $\{H_j\}_{j=1}^N$ an enumeration of the boundary
hypersurfaces of $X$. We assume that each $H_j$ is a smooth embedded
submanifold with corners in $X$, so we can fix a global defining
function $\rho_j$ for that face, i.e.\ $\rho_j \geq 0$, $H_j = \{\rho_j = 0\}$
and $d\rho_j \neq 0$ there.

The space $\calV_b(X)$ of $b$-vector fields on $X$ consists of all
smooth vector fields $V$ which are arbitrary in the interior of $X$
and which lie tangent to all boundary faces, and hence all corners. It
is helpful to describe these in local coordinates. Any point $p \in
\del X$ lies in a corner of codimension $k$, and hence is in the
intersection of $k$ faces $H_{j_1}, \ldots, H_{j_k}$.  There are local
coordinates $(x_1, \ldots, x_k, y)$ based at $p$ with $x_i = \rho_{j_i}$,
$i=1, \ldots, k$, and where $y = (y_1, \ldots, y_{n-k})$ lies in an open
neighbourhood of $0$ in $\RR^{n-k}$.  In terms of these, any $V \in
\calV_b$ can be written as
\[
V = \sum a_{ij}(x,y) x_i \del_{x_j} + \sum b_\ell(x,y)\del_{y_\ell}, \qquad a_{ij}, b_\ell \in \calC^\infty.
\]
In other words, $\calV_b(X)$ is spanned over $\calC^\infty(X)$ by the
local sections $x_i \del_{x_j}$, $\del_{y_\ell}$, $i,j = 1, \ldots, k$,
$\ell = 1, \ldots, n-k$.

Now define the space of conormal functions of order $0$ on $X$
\[
\calA^0(X) = \{u: V_1, \ldots, V_r u \in L^\infty(X)\ \forall\, r\geq 0\ \mathrm{and}\ V_j \in \calV_b\}.
\]
Any $u \in \calA^0$ is smooth in the interior of $X$, and is in some
sense tangentially regular at all boundaries (though it need not have
a well-defined boundary value). Defining regularity using these
$b$-derivatives rather than ordinary derivatives allows functions like
$x^\gamma$, $\Re(\gamma) \geq 0$, and $|\log x|^s$, $s \in \CC$. Next, for any
$k$-tuple $\sigma = (\sigma_1, \ldots, \sigma_k) \in \RR^k$, set $\rho^\sigma = \rho_1^{\sigma_1}\ldots
\rho_k^{\sigma_k}$, and define
\[
\calA^\sigma(X) = \{u = \rho^\sigma v: v \in \calA^0(X)\}.
\]
This is space of conormal functions of (multi)weight $\sigma$. We also write 
$\calA(X) = \cup_{\sigma} \calA^{\sigma}(X)$.

A particularly useful subclass of conormal functions is comprised by
the polyhomogeneous functions.  By definition, $u \in \calA(X)$ is
polyhomogeneous if near any point $p$ on a corner of codimension $k$
in $X$, $u$ has an asymptotic expansion
\[
u \sim \sum a_{\gamma,\ell}(y)x^\gamma (\log x)^\ell, \qquad a_{\gamma,\ell} \in \calC^\infty.
\]
Here $\gamma = (\gamma_1,\ldots, \gamma_k) \in \CC^k$, $\ell = (\ell_1,\ldots, \ell_k)$,
\[
x^\gamma = x_1^{\gamma_1}\ldots x_k^{\gamma_k}, \qquad (\log x)^\ell = |\log x_1|^{\ell_1} \ldots |\log x_k|^{\ell_k},
\]
and where $(\gamma,\ell)$ vary over a discrete set of values in $\CC^k ×
\NN^k$ which has finite intersection with each sector $\cup_{j=1}^k
\{\gamma: \Re (\gamma_j) \leq C_j\} × \NN^k$.

We can also talk about conormality or polyhomogeneity of a function at
an interior $p$-submanifold $Y \subset X$, which by definition is locally
of the form $\{x_{i_1} = \ldots = x_{i_j} = 0, \ y_{r_1} = \ldots = y_{r_s} =
0\}$ (the $p$ means `product', i.e.\ $Y$ is locally of product form in
$X$). We reduce to the situation above by blowing up $X$ along $Y$, so
that $u$ is conormal or polyhomogeneous at $Y$ if it has either of
these properties at the boundary face of $[X;Y]$ which covers $Y$.

Although most of the specific functions we deal with below are
actually polyhomogeneous, we shall not emphasize or need this extra
structure, except in the following weak sense. Fix a boundary face $H$
and suppose that the weight $\sigma_H$ corresponding to this face is $0$.
Let us say that $u \in \calA^{\sigma}_H(X)$ if $u \in \calA^\sigma(X)$ as
before, and that near $H$, $u = u_0 + v$ where $u_0, \rho_H^{-\e}v \in
\calA^\sigma$ and $u_0$ is smooth up to $H$. This is simply a way of
saying that $u$ decomposes into a `leading coefficient' $u_0$, which
is a smooth function on $H$ (conormal at all boundaries of $H$) and a
conormal remainder term $v$ which vanishes to some positive order. If
$\calH$ is a subset of the set of all boundary faces, then
$\calA^\sigma_{\calH}(X)$ consists of functions with this type of
decomposition at each face $H \in \calH$.

Finally, let $Y$ be an interior $p$-submanifold and $\calH$ some
subset of the boundary faces, all elements of which intersect $Y$ at
$\del X$. As before, suppose that $\sigma_H = 0$ for all $H \in \calH$.
Define $\calA^{\sigma}_{\calH}(X;Y)$ to consist of the functions $u$ which can
be decomposed as a sum $u' + u''$ where $u' \in \calA^{\sigma}_\calH(X)$
and $u''$ is supported in a small neighbourhood of $Y$ and
polyhomogeneous on $Y$, and such that this polyhomogeneous singularity
of $u''$ is smoothly extendible across all boundaries of $Y$.

All of these definitions generalize immediately if $u$ is a section of
some smooth vector bundle over $X$.

\subsection*{$\KT$ pseudodifferential operators}
The $\KT$ double space provides the geometric setting for Schwartz
kernels of $\KT$ pseudodifferential operators.
\begin{defi} For any $\mu \in \RR$ and set of weights $\sigma =
  (\sigma_{10},\sigma_{01},0)$ corresponding to the boundary faces $B_{10}$,
  $B_{01}$ and $B_{11}$ of $M^2_{\KT}$, the space
  $\Psi^{\mu,\sigma}_{\KT}(M)$ of conormal $\KT$ pseudodifferential operators
  on $M$ consists of all those operators $A$ on $M$ with the following
  properties:
\begin{itemize}
\item the Schwartz kernel $K_A$ of $A$ is the pushforward (under the
  blowdown $\beta: M^2_{\KT} \to M^2$) of a distribution $\kappa_A$ on
  $M^2_{\KT}$;
\item $\kappa_A$ is a distribution on $M^2_{\KT}$ which is conormal with
  respect to all boundaries, partially homogeneous with respect to the
  front face, and which has a polyhomogeneous singularity of
  pseudodifferential order $\mu$ along the lifted diagonal
  $\mathrm{diag}_{\KT}$, i.e.
\[
\kappa_A \in \calA^\sigma_{B_{11}}(M^2_{\KT},\mathrm{diag}_{\KT}).
\]
\end{itemize}
Slightly more generally, we also define $\Psi^{\mu,\sigma}_{\KT}(M)$ in an
analogous way when $\sigma = (\sigma_{10}, \sigma_{01}, \sigma_{11})$ and 
$\sigma_{11} > 0$ (however, dropping the partial polyhomogeneity at $B_{11}$
and only requiring that $\kappa_A$ is conormal and vanishes to order $\sigma_{11}$
at that face). 
\end{defi} 
The action of $K_A$ on a function $f$ on $M$ requires the choice of a
density $\gamma$ on $M$ against which to integrate, so that $(K_A f)(z) =
\int_M K_A(z,z')f(z')\, \gamma(z')$. It is purely a matter of convention
whether we fix $\gamma$ to be a smooth measure on $\overline{M}$, for
example, or some power of a defining function times a smooth measure;
any two such choices yield equivalent theories, but one does need to
make an adjustment to the index set $\sigma$ below based on this choice.  
We shall follow the convention that $\gamma$ is instead a volume form for some fixed $\KT$
metric, and hence is of the form $x^{-n-d}$ times a smooth measure.
This has the advantage that most of the operators in our later
applications are self-adjoint.

Basic facts about these spaces of pseudodifferential operators include
the composition law
\begin{equation}
\Psi^{\mu,\sigma}_{\KT} \circ \Psi^{\mu',\sigma'}_{\KT} \subset \Psi^{\mu+\mu', \sigma''}_{\KT},\qquad
\sigma'' = (\sigma_{10}, \sigma'_{01},0),
\label{eq:ktcomp}
\end{equation}
which holds provided $\sigma_{01} + \sigma'_{10} > -1$ (this condition is
needed to ensure that the integration defining the composition makes
sense), and the existence of a short exact symbol sequence
\begin{equation}
0 \to \Psi^{\mu-1,\sigma}_{\KT} \to \Psi^{\mu,\sigma}_\KT \to S^\mu({}^{\KT}T^*M) \to 0.
\label{eq:ss}
\end{equation}
The other key facts we need concern the mapping properties of these
operators. To state these we first describe the appropriate function
spaces. Fix any smooth $\KT$ metric $g$ on $M$. This determines the
space $L^2(M;dV_g)$, as well as the basic Hölder space
$\Lambda^{0,\alpha}_{\KT}(M)$, which by definition is the closure of bounded
$\calC^\infty$ functions with respect to the norm
\[
||u||_{0,\alpha} = \sup_{p \in M} |u(p)| + 
\sup_{p \neq q \atop {\mathrm{dist}}_g(p,q) \leq  1} \frac{|u(p)-u(q)|}{\mathrm{dist}_g(p,q)^\alpha}.
\]
Next, for any $s \in \NN$, set
\[
H^s_{\KT}(M) =  \{u: V_1 \ldots V_\ell u \in L^2(M;dV_g): V_j \in \calV_{\KT}\ \forall\, j \leq \ell,\ \ell \leq s\}. 
\]
\[
\Lambda^{s,\alpha}_{\KT}(M) = \{u: V_1 \ldots V_\ell u \in \Lambda^{0,\alpha}_{\KT}: 
V_j \in \calV_{\KT}\ \forall\, j \leq \ell,\ \ell \leq s\}.
\]
and finally, for any defining function $x$ for $\del M$ and $\delta \in \RR$,
\[
x^{\delta} H^s_{\KT}(M) =  \{u = x^{\delta} v: v \in H^s_{\KT}(M)\}, 
\]
\[
x^{\nu}\Lambda^{s,\alpha}_{\KT}(M) = \{u = x^{\nu} v: v \in \Lambda^{s,\alpha}_{\KT}(M)\}.
\]
There is a somewhat loose relationship between certain of these
weighted Sobolev and Hölder spaces.  This is based on the fact that
$x^a \in x^\delta H^s_{\KT}$ (at least locally near $x=0$) if and only if
$a > \delta + (n+d-1)/2$; similarly, $x^a \in x^\nu \Lambda^{s,\alpha}_\KT$ near
$x=0$ if and only if $a \geq \nu$.  Because of this we say that $x^\delta
H^s_{\KT}$ and $x^\nu \Lambda^{s,\alpha}_{\KT}$ are commensurable when $\nu = \delta
+ (n+d-1)/2$.

\begin{prop} Fix $\mu \in \NN$, $\sigma_{10}, \sigma_{01}, \delta \in \RR$ such that $\sigma_{01} + \delta > 0$,
$\sigma_{10} > \delta$. Let $A \in \Psi^{\mu,\sigma}_{\KT}(M)$. Then the maps 
\[
\begin{array}{rcl}
A: x^{\delta} H^{k+\mu}_{\KT}(M) & \longrightarrow & x^{\delta} H^k_{\KT}(M) \\
A: x^{\delta + (n+d-1)/2} \Lambda^{k+\mu,\alpha}_{\KT}(M) & \longrightarrow & x^{\delta +(n+d-1)/2} \Lambda^{k,\alpha}_{\KT}(M)
\end{array}
\]
are bounded.
\label{pr:bddness}
\end{prop}

We sketch a few points in the proof. First note that since these
spaces are defined relative to $\KT$ derivatives, and since
$\calV_{\KT} \circ \Psi^{j}_{\KT} \subset \Psi^{j+1}_{\KT}$, we may immediately
reduce to the case $\mu=0$. Furthermore, conjugating by $x^{-\delta}$ in
the first case and $x^{-(\delta+ (n+d-1)/2)}$ in the second, and observing
that $\tilde{x}/x$ lifts to a conormal function on $M^2_{\KT}$ which
is smooth (and nonvanishing) up to $B_{11}$ reduces us further to the
unweighted case. Finally, the boundedness of $\Psi^0_{\KT}(M)$ on
$L^2(M;dV_g)$ may be deduced via Hörmander's method of using the
symbol calculus to find $B \in \Psi^0_{\KT}$ satisfying $A^* A + B^* B =
C^2 \, \mathrm{Id} + R$ for some $C > 0$ and $R \in \Psi^{-\infty}_{\KT}$ and
then using Cauchy-Schwarz to prove $L^2$ boundedness for operators of
order $-\infty$. The boundedness of $\Psi^0_{\KT}$ on $\Lambda^{0,\alpha}_{\KT}$ is
also deduced in two steps. Decompose $A$ into a sum $A' + A''$ where
the Schwartz kernel of $A'$ vanishes to infinite order at $B_{10}$ and
$B_{01}$ and that of $A''$ is smooth across the diagonal. The
boundedness of $A'$ on Hölder spaces is equivalent to the standard
local boundedness of pseudodifferential operators (of order $0$) on a
neighbourhood in $\RR^n$, cf.\ \cite{Tay96b}.  This argument is
discussed in detail in \cite{Maz91} and \cite{Mazzeo-EMS} for the case
$\KK = \RR$.

\begin{prop} Suppose that $A \in \Psi^{\mu,\sigma}_{\KT}(M)$, where $\sigma =
  (\sigma_{10},\sigma_{01},\sigma_{11})$ indicates conormal order of vanishing at
  each of the three boundary faces $B_{10}$, $B_{01}$, $B_{11}$. If
  $\mu < 0$ and $\sigma_{11} > 0$, and if $\sigma_{01} + \delta > 0$, $\sigma_{10} >
  \delta$, then $A$ is compact on $x^{\delta}L^2$ and on $x^{\delta + (n+d-1)/2}
  \Lambda^{0,\alpha}_{\KT}$.
\label{pr:cpctness}
\end{prop}
This follows directly from the Arzela-Ascoli theorem.

\subsection*{The parametrix construction for fully elliptic operators}
We finally apply the theory of $\KT$ operators outlined above to prove
that under certain hypotheses, the linearized gauged Einstein operator
Laplacian $L^g$ is an isomorphism on certain weighted $L^2$ and Hölder
spaces.  We state the result in slightly greater generality for an
arbitrary generalized Laplacian $P = \nabla^* \nabla + R \in
\Diff^2_{\KT}(M;E)$, associated to a $\KT$ metric $g$. The symmetry of
this operator with respect to the volume form $dV_g$ simplifies some
of the numerology below, but all of these results have direct
analogues for more general fully elliptic $\KT$ operators.

\begin{defi}
  The generalized Laplacian $P^g$ is fully elliptic if its $\KT$
  symbol is invertible as a section of
  $\calC^\infty({}^{\KT}T^*{\overline{M}};\mathrm{Hom}(E))$, and if, in
  addition, its normal operator $N(P^g)$, which is identified via
  Proposition \ref{pr:idno} with the corresponding operator $P^{g_0}$
  on the model hyperbolic space $\KK H^m$, is invertible as an
  unbounded operator on $L^2(\KK H^m; E_0)$.
\label{def:fe}
\end{defi}

Before stating the main theorem of this subsection, let us explore the
relationship of this full ellipticity condition with the indicial root
structure of $P^g$. The indicial roots of $P^g$ are the roots of the
indicial polynomial for $P^g$, and give the rates of vanishing of
formal solutions of this operator. If $\zeta$ is an indicial root, then
there exists some $\phi(y)$ such that $P(x^\zeta \phi(y)) = \calO(x^{\zeta+1})$.
The indicial roots of $P^g$ and of its normal operator $N(P^g)$ are
the same. The indicial roots are arranged symmetrically around
$(n+d-1)/2$ in $\CC$. The complement in $\RR$ of the set of real parts
of all indicial roots of $P$ is a union of open intervals and
half-lines, again symmetric around $(n+d-1)/2$. The significance of
these intervals is as follows. First, if $\delta + (n+d-1)/2$ is the the
real part of some indicial root, hence at the boundary of two
contiguous intervals, then neither of the mappings
\begin{eqnarray}
P: x^{\delta}H^2_{\KT}(M,E;dV_g) & \longrightarrow & x^{\delta}L^2(M,E; dV_g) 
\label{eq:mpPS}\\
P: x^{\delta + (n+d-1)/2} \Lambda^{2,\alpha}_{\KT}(M,E) & \longrightarrow & x^{\delta + (n+d-1)/2} \Lambda^{0,\alpha}_{\KT}(M,E)
\label{eq:mpPH}
\end{eqnarray}
have closed range. This is straightforward to check from basic
definitions. Significantly deeper is the
\begin{theo} Let $P^g\in \Diff^2_{\KT}(M,E)$ be a fully elliptic
  generalized Laplacian. Define $\delta_0$ by the condition that
  $(n+d-1)/2 ± \delta_0$ are the real parts of the indicial roots of $P$
  closest to $(n+d-1)/2$. If $\delta_0 > 0$ and $|\delta| < \delta_0$, then
  (\ref{eq:mpPS}) and (\ref{eq:mpPH}) are both Fredholm, and are
  isomorphisms if and only if the nullspace of $P^g$ on $L^2(M;dV_g)$
  is trivial.
\label{th:fe}
\end{theo}
\begin{proof} The fact that these mappings are Fredholm when $|\delta| <
  \delta_0$ will follow immediately if we can establish the existence of a
  parametrix $G \in \Psi^{-2,\sigma,0}_{\KT}(M;E)$ for $P$ with the property
  that $PG - I = GP - I = Q \in \Psi^{-\infty, \sigma,\infty}_{\KT}(M;E)$. Here for
  convenience we set $\tau = (n+d-1)/2 + \delta_0$ and $\sigma =
  (\sigma_{10},\sigma_{01}) = (\tau,\tau)$. The final index (here $0$ or $\infty$)
  corresponds to $\sigma_{11}$.

  This parametrix is constructed in stages. We first choose an element
  $G_0$ in the small calculus, i.e.\ $G_0 \in \Psi^{-2}_{\KT}(M;E)$, so
  that $PG_0 = I - Q_0$, $Q_0 \in \Psi^{-\infty}_{\KT}$.  This uses only the
  symbol calculus and the symbol ellipticity of $P$, and proceeds
  exactly as in the usual (local) elliptic parametrix construction.

  For the second step we seek a correction term $G_1 \in
  \Psi^{-\infty,\sigma,0}_{\KT}$ chosen so that the remainder term $Q_1 = I -
  P(G_0 + G_1)$ lies in $\Psi^{-\infty,\sigma,1}_{\KT}$, and in particular is
  compact. For this we must solve the normal problem $N(P)N(G_1) =
  N(Q_0) \in \calC^\infty_0(\KK H^m;E)$.  By the second part of the full
  ellipticity hypothesis there is a unique solution to this equation
  in $L^2(\KK H^m,E)$ and (since the right hand side is $\calC^\infty_0$)
  it is a simple matter to check that the solution is conormal at the
  boundaries of the quarter-sphere fibres of the front face. Indeed,
  using the analysis from \cite{Biq00}, we obtain that $N(G_0) \in
  \calA^{\sigma}(S^{n-1}_{++})$, where $\sigma = (\tau,\tau)$ gives the orders of
  conormal vanishing at the two boundaries of the quarter-sphere
  fibres.

  Using $|\delta| < \delta_0$, we see from Proposition \ref{pr:bddness} that
  $G_1$ and $Q_1$ are bounded between these weighted spaces; from
  Proposition \ref{pr:cpctness} we obtain also that $Q_1$ is compact.
  This already shows that $P^g$ is Fredholm. However, it is useful to
  refine this parametrix further.

  Using the composition formula for $\KT$ pseudodifferential
  operators, we see that the iterated compositions of this error term
  with itself vanish to increasingly high order at the front face,
  specifically $Q_1^j \in \Psi^{-\infty,\sigma,j}_{\KT}$. We can therefore take
  an asymptotic sum of the series $R \sim \sum_{j=1}^\infty Q_1^j$ as an
  element of $\Psi_{\KT}^{-\infty,\sigma,1}$.  Now multiply $P(G_0 + G_1) = I -
  Q_1$ on the right by $I+R$. We see that $G' = (G_0 + G_1)(I+R)$
  satisfies $PG' = I - Q'$ where $Q' \in \Psi^{-\infty,\sigma,\infty}_{\KT}$.

  This error term lies in the very residual space of smoothing
  operators with Schwartz kernels which are conormal on $M^2$. One
  consequence is that elements of the nullspace of $L$ in either of
  these function spaces are conormal and vanish like $x^\tau$.
  Furthermore, since these very residual operators form a semi-ideal
  (on one of these weighted $L^2$ spaces, say), a standard argument
  (cf.\ \cite{Maz91}) shows that the true generalized inverse $G$ of
  any of the maps between weighted spaces, which a priori is only
  defined as a bounded operator, is actually an element of
  $\Psi^{-2,\sigma,0}_{\KT}$; the error term $PG - I = GP - I = Q$ is the
  projector onto the nullspace and is still very residual. If the
  nullspace is trivial, then $Q=0$ and hence $P$ is invertible.
\end{proof}

We remark, but do not prove, that if $\delta$ lies in any of the other
open intervals or half-lines described above, then these maps have
closed range but are not Fredholm since either the kernel or cokernel
is infinite dimensional.

With not much more effort, we can prove that the Schwartz kernel of
$G$ has a polyhomogeneous expansion at all boundary faces of
$M^2_{\KT}$ (as well as a polyhomogeneous expansion along the lifted
diagonal of this space which is smoothly extendible across the front
face). This implies that $G$ maps polyhomogeneous sections to
polyhomogeneous sections, and also shows that any section $\kappa$ which
satisfies $L\kappa = 0$ (even just in a neighbourhood of infinity) must
have a complete polyhomogeneous expansion there.

\subsection*{The resolvent family}
The invertibility of the linearized gauged Einstein operator $L^g$ on
weighted Hölder spaces, which is a direct consequence of Theorem
\ref{th:fe} and Proposition \ref{pr:bddness}, is the key ingredient in
the deformation theory of A$\KK$H Einstein metrics. For the analogous
result on products of A$\KK$H spaces, we shall use a spectral
synthesis formula for the inverse of this operator which expresses it
in terms of the resolvent families of the corresponding operators on
each factor. In preparation for this, we now recall the basic theory
of these resolvent families in terms of the $\KT$ calculus and prove
some estimates on their Schwartz kernels which are uniform in the
spectral parameter.

First let us recall some results from \cite[§~I.4]{Biq00}. Define $\delta_0 = \delta_0^\KK$
for the operator $L^g$ as in the statement of Theorem \ref{th:fe}; then
\begin{equation}
\label{eq:2}
\begin{array}{rlcl}
  \delta_0 &= (n+d-1)/2 & \text{for }\bK&=\RR \text{ or }\CC,\\
  \delta_0 &> (n+d-1)/2 &   \text{for }\bK&=\HH \text{ or }\OO.
\end{array}
\end{equation}
This means that the interval of weights $\delta$ for which $L^g$ is
Fredholm is exactly $(0,n+d-1)$ in the real and complex cases, and is
larger in the quaternionic and octonionic cases.

By definition, the resolvent of $L^g$ is the family of $L^2$ bounded
operators $(L^g - \lambda)^{-1}$, which exists precisely when $\lambda \notin
\mathrm{spec}\,(L^g)$. We wish to recognize these operators as
elements of $\Psi^{-2,*}_{\KT}(M)$, depending holomorphically on $\lambda$ in
an appropriate sense. This will follow from Theorem \ref{th:fe}, which
in turn requires the
\begin{lemm}
  The operator $L^g - \lambda$ is fully elliptic (as a $\KT$ operator) if
  and only if $\lambda \notin [\delta_0^2,\infty)$.
\label{le:feo}
\end{lemm}
\begin{proof} First note that 
\[
{}^{\KT}\sigma_2(L^g-\lambda) = {}^{\KT}\sigma_2(\nabla^* \nabla) = |\zeta|^2\, \mathrm{Id},
\]
which is obviously invertible. In addition, $N(L^g - \lambda) = L^{g_0} -
\lambda$, so we conclude that $L^g-\lambda$ is fully elliptic if and only if $\lambda
\notin \mathrm{spec}\,(L^{g_0})$.

The indicial operator of $L^{g_0}$ is a second order matrix-valued
ordinary differential operator, and the indicial roots correspond to
solutions of the form $x^\zeta \kappa_0$, where $\kappa_0$ is a constant
symmetric two-tensor. By reducing to the various irreducible
components in $\Sym^2$, we obtain them as the roots of a
finite number of quadratic polynomials $\zeta^2 - (n+d-1) \zeta + \alpha$, where
$\alpha$ is a constant depending on dimension and the irreducible
component of the decomposition. The roots from any one of these
polynomials are $(n+d-1)/2 ± \frac12 \sqrt{(n+d-1)^2 - 4\alpha}$. For some
$\alpha_0$ we obtain the roots with real part closest to $(n+d-1)/2$, that
is $(n+d-1)/2±\delta_0$. So we see that each $\alpha \leq \alpha_0$. Now, the
indicial roots of $L^{g_0} - \lambda$ are the roots of the various
polynomials $\zeta^2 - (n+d-1)\zeta + \alpha + \lambda$, hence are equal to 
\[
(n+d-1)/2 ± \frac12\sqrt{(n+d-1)^2 - 4\alpha - 4\lambda}.
\]
Define
\[
\delta_0(\lambda) = \min_\alpha \frac12 \Re \sqrt{(n+d-1)^2 - 4\alpha - 4\lambda} \ 
= \ \Re \, \sqrt{\delta_0^2 - \lambda}.
\]
By Theorem \ref{th:fe}, $L^{g_0}-\lambda$ is at least Fredholm on $L^2$ provided
\[ 
\delta_0(\lambda) > 0, 
\]
or equivalently, if $\lambda \notin [\delta_0^2,\infty)$. 

So far we have proved that if $\lambda$ is outside this half-line, then
$L^{g_0}-\lambda$ is at least Fredholm.  This shows that
$\mathrm{spec}\,(L^{g_0})$ is the union of $[\delta_0^2,\infty)$ and finite
point spectrum of multiplicity in the half-line $(-\infty, \delta_0^2)$.
However, this point spectrum must be empty, since otherwise, if 
$L^{g_0} \phi = \widehat{\lambda} \phi$, for some $\widehat{\lambda} < \delta_0^2$
and with $\phi \in L^2$, then the subspace spanned by all translates of $\phi$ by 
isometries of $\KK H^m$ would be infinite dimensional, contradicting the fact 
that $L^{g_0} - \widehat{\lambda}$ is Fredholm.  This finishes the proof.
\end{proof}

The same reasoning leads to the 
\begin{theo} Let $L^g$ be the linearized gauged Einstein operator on
  the manifold $M$ with $\KT$ metric $g$. Then $L^g - \lambda$ is Fredholm
  if and only if $\lambda \notin [\delta_0^2,\infty)$. More precisely,
\[
\mathrm{spec}\,(L^g) = [\delta_0^2,\infty) \cup \{\lambda_i\}_{i=1}^N,
\]
where $\lambda_i$ lies in $(-\infty,\delta_0^2)$ and is an $L^2$ eigenvalue of
finite multiplicity.
\label{th:resfam1}
\end{theo}
An A$\KK$H Einstein space $(M,g)$ is nondegenerate if and only if $0$
is not in this point spectrum.

We shall need to know slightly more about the dependence of the inverse on $\lambda$.
\begin{prop}
  Fix $\epsilon > 0$ and define $\Omega_\epsilon \subset \CC \setminus [\delta_0^2,\infty)$ to
  consist of the set of all $\lambda$ for which $\delta_0(\lambda) > \epsilon$. Let $\tau_\e
  = (n+d-1)/2 + \e$ and $\sigma_\e = (\tau_\e,\tau_\e,0)$.  Then for each $\e
  > 0$, the resolvent family
\[
\Omega_\e \ni \lambda \mapsto R(\lambda) = (L^g-\lambda)^{-1} \in \Psi^{-2,\sigma_\e}_{\KT}(M;\Sym^2({}^{\KT}T^*M))
\]
is meromorphic in the sense that the Schwartz kernels of these
operators, as elements of a fixed space of distributions, depends
meromorphically on $\lambda$. The poles occur only at each $\lambda_i$ and these
are all simple; the residues are the finite rank orthogonal
projections onto the corresponding eigenspaces.
\label{pr:resfam2}
\end{prop}
The proof is based on the fact that the model resolvent
$(L^{g_0}-\lambda)^{-1}$ is itself holomorphic, which can be checked by
direct ODE analysis, and the analytic Fredholm theorem. This is the
direct generalization of \cite{MazMel87} and \cite{EpsMelMen91}, cf.\
also \cite{Gui05}.

We conclude this section by proving uniform estimates for the
off-diagonal Schwartz kernel of this resolvent when $\lambda = i\mu$ lies on
the imaginary axis.
\begin{prop} Let $(M,g)$ and $L^g$ be as above. When $\mu \in \RR$, the
  indicial root of $L^g - i\mu$ which has the smallest strictly
  positive real part is equal to
\[
\zeta(\mu) = \frac{n+d-1}{2} + \sqrt{\delta_0^2 - i\mu}. 
\] 
In particular $\gamma(\mu) := \Re\,\zeta(\mu) \geq n+d-1$ with equality if
and only if $\mu = 0$.  Furthermore,
\[
\zeta(\mu) \sim \sqrt{|\mu|/2}\,(1 ± i) \qquad \mathrm{as} \quad \mu \to ± \infty.
\]
Let $K(z,z',i\mu)$ be the Schwartz kernel of this operator. Then there
exists $\mu_0 > 0$ such that for any $c > 0$ and $\e \in (0,1)$, if
$z,z' \in M$ satisfy $d(z,z') \geq c > 0$ and $|\mu| \geq \mu_0$, we have
\[
|K(z,z',i\mu)| \leq C \mu^{-1} x^{(1-\e)\Re\, \zeta(\mu)}
\]
where the constant $C$ is independent of $\mu$. A similar estimate
holds for all $b$-derivatives of $K$.
\label{pr:expbounds}
\end{prop}
\begin{proof}
  The calculation of the indicial root and the statements about its
  asymptotics are straightforward, based on the remarks in the proof
  of Lemma \ref{le:feo}.  As for the main assertion, first suppose
  that $f \in \calC^\infty_0(M)$, $|f| \leq 1$, and define $u = u_\mu =
  R(i\mu)f$. We claim that for any $\e > 0$, there exists a constant
  $C_\e > 0$ such that for $|\mu| \geq \mu_0$,
  \[
  |u| \leq C \mu^{-1} x^{(1-\e)\gamma(\mu)}.
  \]
  In particular, $C$ is independent of $\mu$ (and in addition, $\sup
  |u|$ depends linearly on $\sup |f|$).

  To prove this, first recall that since $R(i\mu) \in
  \Psi^{-2,\gamma(\mu),\gamma(\mu)}_{\KT}(M)$, it is immediate that $|u| \leq A_\mu
  x^{\gamma(\mu)}$; the issue is to prove the uniformity in $\mu$ of the
  constant $A$. Suppose this fails, i.e.\ suppose there exists a
  sequence of $\mu$ tending to infinity so that
  \[
  \sup_{z \in M} \mu |u(z)| x^{-(1-\e)\gamma(\mu)} = A_\mu \to \infty.
  \]
  This supremum is attained at a point $q_\mu \in M$, so if we define
  \[
w(z) = \left(\mu x(q_\mu)^{-(1-\e)\gamma(\mu)}/A_\mu\right) u(z),
  \]
then
\begin{equation}
|w(z)| \leq (x/x(q_\mu))^{(1-\e)\gamma(\mu)},
\label{bdw}
\end{equation}
with equality at $z = q_\mu$ and 
\begin{equation}
(\mu^{-1}L - i) w(z) = \left(x(q_\mu)^{-(1-\e)\gamma(\mu)}/A_\mu\right)f(z).
\label{eqw}
\end{equation}
We shall consider various cases depending on whether or not $q_\mu$
remains in a compact set of $M$.

Suppose first that $x(q_\mu) \geq c > 0$. Then the right hand side of
(\ref{eqw}) tends to zero uniformly.  Let $B$ be a geodesic ball of
radius $1$ centered at $q_\mu$, fix a trivialization of the bundle over
$B$ and suppose that $z$ are Riemann normal coordinates in this ball.
Set $\xi = \sqrt{\mu}z$, so that $\xi$ lies in a ball of radius
$\sqrt{\mu}$ in $\RR^{n+1}$. In terms of these coordinates,
\[
\mu^{-1}L = -\sum_{j=1}^{n+1}\frac{\del^2\,}{\del \xi_j^2} + \calO(\mu^{-1/2}).
\]
The remainder term is a second order operator with coefficients which
converge to $0$ uniformly on compact sets in these expanding balls in
$\RR^{n+1}$. Using standard local elliptic theory, we can take a limit
in (\ref{eqw}), and obtain a function $w_\infty$ defined on the entire
Euclidean space, such that
\[
(\Delta_\xi - i) w_\infty = 0, \qquad |w_\infty| \leq 1,\ |w_\infty(0)| = 1.
\]
However, no such function exists. To see this, take Fourier transform
(for $w$ as an element of ${\mathcal S}'(\RR^{n+1})$) and use that the
full symbol $|\xi|^2 - i$ is nowhere vanishing.  Hence this case cannot
occur.

Now suppose that $x(q_\mu) \to 0$. Pass to a subsequence so that $q_\mu \to
\bar{q} \in \del M$, and then apply a sequence of parabolic dilations
$D_\mu$ based at appropriate points converging to $\bar{q}$ and with
strength $x(q_\mu))^{-1}$ so that $D_\mu(q_\mu) = (1,0)$ in a fixed
coordinate system $\tilde{z}$.  We now proceed much as before. Let $B$
be a unit geodesic ball centered at $(1,0)$, and define $\xi =
\sqrt{\mu}\tilde{z}$. The sequence of operators $\mu^{-1}D_\mu^*L$
converge to $\Delta_\xi$ as before.  The bound on $\tilde{w} = D_\mu^* w$
now takes the form
\[
|\tilde{w}| \leq e^{(1-\e)\gamma(\mu)/\sqrt{\mu}}
\]
with equality at the origin; here $t = \log (x/x(q_\mu))$, and we can
assume this is the first coordinate $\xi_1$ in the $\xi$ system. We
again pass to a limit. The limiting function $\tilde{w}_\infty$ satisfies
$|\tilde{w}_\infty(0)| = 1$,
\[
(\Delta_{\xi} - i)\tilde{w}_\infty = 0, \qquad |\tilde{w}_\infty| \leq e^{(1-\e)t/\sqrt{2}}.
\]
To analyze whether this is possible, note that this exponential bound
on $\tilde{w}_\infty$ implies that its Fourier transform is well defined
as an element of ${\mathcal S}'$ on the subspace $\{\xi \in \CC^{n+1}:
\Im\, \xi_1 = (1-\e)/\sqrt{2}, \xi_j \in \RR, j > 1\}$.  The symbol $\xi
\cdot \xi - i$ is again invertible here, which precludes the existence of
this limit; hence this case is also impossible. This proves that the
function $u = R(i\mu)f$ satisfies the stated bound uniformly in $\mu$.

An essentially identical argument proves that a similar bound holds
regardless of the location of the support of $f$. In other words,
suppose that $\mathrm{supp}\,(f) \subset B_1(p_\mu)$ and $\sup |f| \leq 1$.
Then for the $L^2$ solution to $(L-i\mu)u = f$, we have $|u(q)| \leq C
\mu^{-1} \exp (-(1-\e)d(q,p_\mu))$.  The only modification needed is
that if there a sequence $u_\mu$ for which the constant increases
without bound, and if the center $q_\mu$ of the support of $f_\mu$ tends
to infinity, then we parabolically rescale so as to obtain a sequence
of problems $(L_\mu - i\mu)\tilde{u}_\mu = \tilde{f}_\mu$, where the
rescaled operators $L_\mu$ converge to the limiting model operator for
the $\KT$ structure. The validity of the bound in this case follows by
what we have done above.

We have now proved that $K(z,z',i\mu)$ decays like
$e^{-(1-\e)d(z,z')}\mu^{-1}$ for $d(z,z') \geq c > 0$ in a weak sense.
More precisely, let $\calU$ be any neighbourhood in $M × M$ with
compact closure which does not intersect the diagonal; then if $\sigma >
2n+2$, the $H^{-\sigma}$ norm of the restriction of Schwartz kernel to
$\calU$ satisfies this bound. Using that $(L_z - i\mu)K = (L_{z'} -
i\mu)K = 0$ away from the diagonal, we can estimate any ${\mathcal
  C}^k$ norm of $K$ in $\calU$ at the cost of introducing an extra
factor $\mu^{k + \sigma}$ into the estimate. This in turn may be absorbed
into the exponential by decreasing the factor $\e$ slightly. This
completes the proof of the $\calC^0$ bound, and indeed also of bounds
with respect to any $\calC^k$ norm in the interior. In fact, it gives
slightly more, namely that this bound holds even after applying any
sequence of $\KT$ vector fields to $K$ on the left and right; this is
because $\KT$ derivatives are controlled by powers of $L$, which as
above are equivalent to powers of $\mu$.

To finish, we also need to check the conormal bounds, i.e.\ that the
same estimates remain true if we apply any sequence of $b$ vector
fields to $K$ on the left and right. For this we point out the
following facts: first, since $K(z',z,i\mu) = K(z,z',-i\mu)^*$, we need
only the case where all $b$ derivatives are applied to the left -- $z$ -- 
factor; next, $K$ itself is conormal, so these $b$-derivatives behave
well locally uniformly in $\mu$, i.e.\ it is only the large $\mu$
behaviour that might be problematic; finally, we can repeat the same
proof as for $K$ itself, using at the final step to convert the weak
bounds to strong ones that if $V$ is any $b$-vector field, then
$[L,V]$ is a $\KT$-operator of order $2$, hence is bounded by
multiplication by $\mu$.
\end{proof}

\begin{rema}
  We have stated the results on the resolvent family for the
  linearized Einstein operator, but the results remain true for any
  geometric Laplacian, provided we choose $\delta_0$ as in Theorem
  \ref{th:fe}. For example, on an asymptotically quaternionic
  hyperbolic space, Theorem \ref{th:resfam1} gives the spectrum of the
  Hodge Laplacian acting on differential forms (except when the
  degree equals half the dimension, then there is a zero eigenvalue of
  infinite multiplicity).
\end{rema}

\section{Einstein deformation theory}\label{sec:einst-deform-theory}
We now present some basic facts about the (Bianchi gauged) Einstein
operator and its linearization.  Using results from the last section,
we review how this yields the deformation theory for the rank one
hyperbolic spaces in the class of A$\bK$H Einstein spaces. This is
contained in \cite{Biq00} for all $\KK$, see also \cite{GraLee91} and
\cite{Lee06} for the result when $\KK = \RR$, so the only novelty here
is showing how this follows immediately through the use of the $\KT$
calculus. These same arguments are used again in §~\ref{sec:sef} for a
coupled generalization of these same equations, and in the product
case in §~\ref{sec:gedt}.

\subsection{The Einstein equation and the Bianchi gauge}
\label{ssec:einst-equat-bianchi}
The Einstein equation $\Ric^g + \lambda g = 0$ is not elliptic because of
its diffeomorphism invariance. Amongst many viable gauge choices, the
Bianchi gauge introduced in \cite{Biq00} is particularly convenient.
Define the map from symmetric $2$-tensors to $1$-forms, relative to
the fixed background metric $g$,
\begin{equation}
k \longmapsto B^g(k) = \delta^g k + \frac12 \, d \, \tr^g \, k.
\label{eq:bg}
\end{equation}
Note that $B^g(g) = 0$, so the subspace of metrics $\tilde{g}$ near to
$g$ which are in Bianchi gauge (i.e.\ so that $B^g(\tilde{g}) = 0$),
is identified with the set of tensors $k$ near $0$ such that $B^g(k) =
0$. The system
\begin{equation}
\Ric(\tg) + \lambda\tg = 0, \qquad B^g(\tg) = 0,
\label{eq:system}
\end{equation}
which is elliptic in the sense of Agmon-Douglis-Nirenberg, can be
rolled up into the single elliptic equation
\begin{equation}
N^g(k) := \Ric(g+k) + \lambda (g+k) + (\delta^{g+k})^*B^g(k) = 0.
\label{eq:main}
\end{equation}
As proved in \cite[chapter 1]{Biq00}, 
\begin{prop}\label{prop:einst-equat-bianchi}
Suppose that $N^g(k) = 0$, and in addition that $|B^g(k)| \to 0$ and the Ricci curvature of $g+k$ is 
nonnegative and strictly negative somewhere. Then $g+k$ satisfies (\ref{eq:system}), i.e.\ is
Einstein and in Bianchi gauge.
\end{prop}
The proof follows from the Weitzenböck formula
\[
B^{g+k}N^g(k) = \delta^{g+k}(\delta^{g+k})^*B^g(k) = \big( (\nabla^{g+k})^*\nabla^{g+k} - \Ric^{g+k} \big) B^g(k)
\]
and the Bochner technique.

{}From the same Weitzenböck formula, the converse follows almost
immediately, that is, any Einstein metric $\tg$ close to $g$ can be
put in the Bianchi gauge to satisfy the system (\ref{eq:system}). More
precisely, let $\mathrm{Diff}$ denote the set of all diffeomorphisms
on the A$\bK$H space $M$ which are close to the identity and
exponentials of vector fields $X \in x^{\nu}\Lambda^{3,\alpha}_{\KT}(M)$, $\calM$
the set of all metrics $\tilde{g} = g+k$ with $k\in
x^{\nu}\Lambda^{2,\alpha}_{\KT}(M)$, and $\calS\subset\calM$ the set of metrics $\tg$
which satisfy $B^g(\tg) =0$. Then one has the following slice
statement \cite[chapter 1]{Biq00}:
\begin{prop}
If $\Ric^g<0$, then the natural map
\[
\mathrm{Diff} × \calS \longrightarrow \calM, \quad (\tilde{g},\Phi) \longmapsto \Phi^*(\tilde{g})
\]
is a local homeomorphism.
\end{prop}

An advantage of this gauge is that the linearization takes the simple form
\begin{equation}
L^g \kappa := \left. 2 DN^g \right|_0(\kappa) = \nabla^* \nabla  \ \kappa - 
2 \overset{\circ}{R} \kappa \ +  \Ric \circ \kappa + \kappa \circ \Ric  +
2 \, \lambda \, \kappa; 
\label{eq:lin}
\end{equation}
here 
\[
(\overset{\circ}{R} \kappa)_{ij} = R_{ipjq}\, \kappa^{pq}, \qquad (\Ric \circ 
\kappa)_{ij} = 
\Ric_i^{\ p}\, \kappa_{pj}, \qquad (\kappa \circ \Ric^g)_{ij} = \kappa_i^{\ p} \,\Ric_{pj},
\]
and all curvatures and covariant derivatives are computed relative to
$g$.  Note in particular that if $\Ric^g = -\lambda g$, then
\begin{equation}
L^g = \nabla^* \nabla - 2 \overset{\circ}{R}.
\label{eq:lgeo2}
\end{equation}

\subsection{Deformation theory for A$\bK$H Einstein spaces}
\label{ssec:deform-theory-abkh}
We now review the basic deformation theory for A$\bK$H Einstein
spaces, proved originally in \cite{GraLee91} when $\bK = \RR$, and in
\cite{Biq00} in the other two cases. We do not discuss the more subtle
aspects of this deformation theory, but restrict attention to the
simpler case of perturbations of {\it nondegenerate} A$\bK$H Einstein
metrics.

\begin{defi} 
  An A$\bK$H Einstein metric $g$ is said to be nondegenerate if the
  $L^2$ nullspace of the linearized Bianchi-gauged Einstein operator
  $L^g$ is trivial.
\label{de:ndg}
\end{defi}
Recall that $g$ is nondegenerate if for any $\nu > -\delta_0^\KK$, the nullspace of 
$L^g$ on $x^{\nu + (n+d-1)/2}\Lambda^{2,\alpha}_{\KT}(M,S^2(T^*M))$ is trivial. This 
follows from the regularity theorem stating that if $L^g \kappa = 0$ and 
$|\kappa| \leq C x^{\nu'}$ for some $\nu' > (n+d-1)/2 - \delta_0$, then 
$\kappa \in \calA^{(n+d-1)/2 + \delta_0}_{\phg}$, and in particular $\kappa \in L^2$.

The significance of this nondegeneracy condition is contained in the
\begin{prop}[\cite{GraLee91}, \cite{Biq00}, \cite{Lee06}]
  \label{prop:deform-theory-abkh}
  Let $g$ be a nondegenerate A$\bK$H Einstein metric with $\calC^{\infty}$
  conformal infinity $\frakc$. Then every $\calC^\infty$ conformal
  infinity datum $\frakc'$ sufficiently close to $\frakc$ in the
  $\calC^{2,\alpha}$ topology is the conformal infinity of an A$\bK$H
  Einstein metric $g'$ such that $g'-g_{\frakc'} \in x^\nu \Lambda^{2,\alpha}_{\KT}$ for
  some $\nu > 0$. (Here $g_{\frakc'}$ is the A$\KK$H metric with conformal infinity
data from Definition 2 in §2.) This metric $g'$ is unique amongst A$\KK$H
  metrics with the specified conformal infinity and such that $||g'-g_{\frakc'}||_{2,\alpha,\nu}$ is small.
\end{prop}

\begin{proof} 
  First define an extension operator which associates to the conformal
  infinity $\frakc'$ an A$\bK$H metric $g_{\frakc'}$. If $\frakc =
  ([\gamma],\eta)$, then we choose a product decomposition $(0,\epsilon) × X$ of a
  collar neighbourhood of $\partial M$ and a radial coordinate $x$ so that
  $g$ has the form $\frac{dx^2+ \gamma}{x^2}+\frac{\eta^2}{x^4}+k$, with
  $k\in x \Lambda^{2,\alpha}_{\KT}$. (We can take the weight $\nu$ to equal $1$ when
  $\frakc$ is smooth, but could also use any smaller positive value.) 
  Fixing a cutoff function $\chi(x)$ which equals $1$ for $x \leq
  \epsilon/3$ and vanishes for $x\geq 2\epsilon/3$, then for any conformal infinity
  $\frakc' = ([\gamma'],\eta')$, set
\[
g_{\frakc'} = (1-\chi(x))g + \chi(x)\left(\frac{dx^2+\gamma'/4}{x^2}+\frac{(\eta')^2}{4x^4}+k \right).
\]
(In the quaternionic case, the function spaces $x\Lambda^{2,\alpha}_{\KT}$
vary with $\frakc'$ so one has also to choose $k$ varying continuously
with $\frakc'$). It is straightforward that
\[
\|N^g(g_{\frakc'})\|_{0,\alpha,\nu} \leq C \left(\|\gamma' - \gamma\|_{2,\alpha} + 
\|\eta' - \eta\|_{2,\alpha}\right).
\]

Writing $N^g(g_{\frakc'}) = f_{\frakc'}$, then Taylor expansion gives
\[
N^g(g_{\frakc'} + k) = f_{\frakc'} + L^{\frakc'} k + Q(\frakc',k),
\]
where the second term on the right is the linearized Bianchi-gauged
Einstein operator at $g_{\frakc'}$.  The nondegeneracy of $L^g$
implies that $L^{\frakc'}$ is also invertible if $\frakc'$ is close
enough to $\frakc$; we denote its inverse by $G^{\frakc'}$.  When
$\frakc'$ is $\calC^\infty$, this operator is a $\KT$ pseudodifferential
operator of order $-2$, and for $0 < \nu \leq 1$,
\[
G^{\frakc'}: x^\nu \Lambda^{0,\alpha}_{\KT}(M, S^2(T^*M)) \longrightarrow 
x^\nu \Lambda^{2,\alpha}_{\KT}(M, S^2(T^*M))
\]
is bounded. Furthermore, the norm of this operator is bounded
independently of $\frakc'$ in a neighbourhood of $\frakc$.  When $\bK
= \bH$, the function spaces vary with $\eta'$ (i.e.\ with the
distribution $\calD'$).

Now write the equation to be solved as 
\[
k = -G^{\frakc'}\left(f_{\frakc'} + Q(\frakc',k)\right).
\]
The right side defines a contraction mapping when $\frakc'$ is
sufficiently close to $\frakc$, and from this we immediately obtain a
unique solution $k$.  The metric $g_{\frakc'} + k$ is an A$\bK$H
metric which solves the gauged Einstein equation. By Proposition
\ref{prop:einst-equat-bianchi}, it is in fact an A$\bK$H Einstein
metric in Bianchi gauge.
\end{proof}

\begin{rema}
  A precise statement about the regularity of $g$ and $g'$ near $\del
  X$ has been omitted, and indeed this is a subtle issue.  There is a
  substantial difference between understanding the dependence of the
  asymptotic regularity for an arbitrary A$\KK$H Einstein metric on
  that of its conformal infinity data, and the same question for such
  a metric obtained by perturbation from one which is a priori known
  to be polyhomogeneous.  The reason is that in the former case one
  needs to deal explicitly with the gauge conditions, while in the
  latter, the gauge choice is part of the setup and the perturbation
  term automatically satisfies a $\KT$ elliptic equation. Here is a
  summary of what is known.

For the nonperturbative case, when $\KK = \RR$ and the conformal
infinity data $\frakc$ of the A$\RR$H Einstein metric $g$ is
$\calC^\infty$, then $g$ is polyhomogeneous, and in fact, in even
dimensions has a smooth conformal compactification, \cite{ChrDelLeeSki05},
\cite{Helliwell}.  The corresponding result has not been proved in the
other cases (except when $\KK = \CC$ and $g$ is Kähler-Einstein,
\cite{LeeMel82}), but is surely true by essentially the same method
as in \cite{ChrDelLeeSki05}.

As for the perturbative case, a simple adaptation of the argument in
\cite{Maz91b}, which depends only on the commutation properties
of $\Psi_{\KT}^*$ with $b$ vector fields on $X$, proves that if the
perturbed conformal infinity data $\frakc'$ is $\calC^\infty$, then the
solution $k$, and hence the metric $g'$, constructed above is
polyhomogeneous.  If $\KK = \RR$ and $\frakc'$ is only $C^{2,\alpha}$,
then \cite{Lee06} proves that if $g$ is sufficiently regular, then
$g'$ has a $\calC^{2,\alpha}$ conformal compactification, but this is not
a sharp statement since the function space $\calC^{2,\alpha}(\del X)$ is
not well adapted to these types of degenerate problems.

These issues will not be emphasized here, and we shall tacitly assume
the polyhomogeneous regularity of A$\KK$H metrics with $\calC^\infty$
conformal infinities. To ease the reader's conscience, however, since
we will not supply the full proof of that fact, note that this issue
only arises in §~\ref{sec:gener-lapl-near}, and one can easily adapt
the arguments to accommodate metrics with lower regularity, as we
discuss briefly there.
\label{re:reg}
\end{rema}

\section{Asymptotically product hyperbolic metrics and their conformal infinities}
\label{sec:pbs}

Let $M_i^{n_i+1}=G_i/K_i$ be an A$\KK_i$H space, $i=1,2$. The boundary
at infinity, $S^{n_i}=K_i/H_i$, is equipped with the standard
$\bK_i$-contact distribution $\cD_i$, which has a conformal $H_i$
structure inducing a compatible conformal class $[\gamma_i]$. The product
hyperbolic space $M = M_1 × M_2$ is a (reducible) rank two symmetric
space with Furstenberg boundary $X = S^{n_1} × S^{n_2}$.  Let $g_i$ be
the standard metric on each factor, so that $\Ric^{g_i} + \lambda_i g_i =
0$. Then
\[
g = \lambda_1 g_1 + \lambda_2 g_2 
\]
is Einstein with 
\[
\Ric^g + g = 0.
\]

We begin this section by describing a class of boundary structures on
$X$, called $(G_1× G_2)$-conformal structures, which constitute the
conformal infinity data for the class of Einstein metrics we
eventually construct.  The problem of extending one of these boundary
structures to a metric on $M$ which is asymptotically Einstein in an
appropriately strong sense is far from immediate. The `obvious'
extension has Einstein tensor vanishing in some sector near $X$, but
not uniformly near infinity. The main goal of this section is to
explore the geometry of asymptotically product hyperbolic metrics in
order to find the correct compatibility conditions for metrics which
are asymptotically Einstein in this stronger sense. Their construction
is carried out in the next section.

\subsection{$(G_1 × G_2)$-conformal structures}
\label{ssec:g_1-g_2-conformal}
Using the notation above, we make the 
\begin{defi}\label{defi:prod-conf}
  A \emph{$(G_1 × G_2)$-conformal structure} on $X$   consists of a pair
  of distributions, each equipped with a conformal class of metrics,
$(\cD_i,[\gamma_i])$, $i=1,2$, such that
\begin{itemize}
\item[(1)] the distributions $\cF_i=\cD_i+[\cD_i,\cD_i]$ are integrable;
\item[(2)] $\cF_1 \oplus \cF_2 = TX$;
\item[(3)] the pair $(\cD_i, \cF_i/\cD_i)$ with induced bracket 
\[ 
[\ , \, ]: \cD_i × \cD_i \to \cF_i/ \cD_i
\] 
is isomorphic to the graded $\bK_i$-Heisenberg algebra;
\item[(4)] $\cD_i$ is equipped with a conformal $H_i$ structure,
  compatible with the bracket, and inducing the conformal metric $[\gamma_i]$; 
\item[(5)] $[\cD_1, \cD_2] \subseteq \cD_1+\cD_2$.
\end{itemize}
\end{defi}

The basic example, of course, is a product structure on $X = X_1 ×
X_2$: here $\cD_1 \subset \cF_1 = TX_1 \oplus \{0\}$ and $\cD_2 \subset \cF_2 =
\{0\} \oplus TX_2$, and each $(\cD_i,[\gamma_i])$ is the pullback of a
$G_i$-conformal structure from $X_i$. Our main focus in this paper is
with perturbations of these product structures. As we now indicate,
there is substantial rigidity in the deformation theory, and nearby
structures retain many vestiges of the product case.

\begin{lemm}\label{lemm:small-def}
For $i=1,2$, let $X_i$ be a compact simply-connected manifold with $G_i$-conformal structure $(\cD_i,[\gamma_i])$. 
Then any small deformation of the product $(G_1 × G_2)$-conformal structure on $X=X_1 × X_2$ 
has the following properties:
\begin{enumerate}
\item the pairs of distributions $(\cD_1,\cD_2)$ and $(\cF_1,\cF_2)$
  remain of product type;
\item if $\bK_i = \RR$ or $\CC$, then for that $i$, $\cD_i$ remains
  fixed (up to a global diffeomorphism), but the deformation of
  $[\gamma_i]$ may depend on both factors in $X_1 × X_2$;
\item if $\bK_i = \HH$, then the distribution $\cD_i$ varies amongst
  quaternionic contact structures, but since the conformal class
  $[\gamma_i]$ is determined by $\cD_i$, any deformation of $[\gamma_i]$
  depends only on the factor $X_i$;
\item if $\bK_i = \OO$, then both $\cD_i$ and $[\gamma_i]$ remain fixed in
  the deformation (and in fact $X_i = S^{n_i}$ and the structure is
  standard).
\end{enumerate}
In particular, the distribution $\cD_i$ can change (modulo
diffeomorphisms) only in the quaternionic case, and the conformal
metric $[\gamma_i]$ may depend on both factors $X_1$ and $X_2$ only if
$X_i$ is real or complex.
\end{lemm}

\begin{proof}
  At the initial product structure, the leaves of the foliation
  corresponding to $\cF_i$ are just the parallel copies of $X_i$.
  After a small deformation, the leaves are covering spaces for $X_i$,
  and since both $X_1$ and $X_2$ are simply-connected, these leaves
  must remain diffeomorphic to $X_i$. In particular, the perturbed
  distributions $\cF_i$ are still equal to the tangent bundles of the
  respective factors.

Now observe that as an immediate consequence of conditions (3) and (5), 
\[
[\cF_1,\cD_2] = [\cD_1 + [\cD_1,\cD_1],\cD_2] \subseteq \cF_1+\cD_2
\]
(by the Jacobi identity), and hence $\cD_2$ is invariant along the
leaves of the foliation corresponding to $\cF_1$; similarly $\cD_1$ is
invariant along the leaves of the foliation for $\cF_2$.

When $\bK_i = \RR$, $\cD_i$ remains equal to the tangent bundle
$TX_i$, while if $\bK_i = \CC$, then by Darboux's lemma, we may still
assume that $\cD_i$ remains fixed in the deformation.

In other words, in this deformation theory, we may as well assume that
the distributions $\cD_i$ and $\cF_i$ remain of product type. The
remaining assertions follow directly from this.
\end{proof}

When $X_1$ or $X_2$ are not simply connected, we shall impose this
product structure as a separate hypothesis:
\begin{defi}\label{defi:global-int} 
  A deformation of a product $(G_1 × G_2)$  -conformal structure on
$X=X_1 × X_2$ is called \emph{globally integrable} if (modulo
diffeomorphism) the foliations $\cF_1$ and $\cF_2$ remain the tangent
spaces of the two factors of $X=X_1 × X_2$.
\end{defi}

It is possible to define $(G_1 × G_2)$-conformal structures on any
closed manifold $X$ of the appropriate dimension. Looking ahead to the
main goals of this paper, one could then try to extend this to an
asymptotically Einstein metric on some manifold $M$ with two boundary
hypersurfaces $F_1$ and $F_2$ and $X$ as its corner of codimension
$2$. However, we have already noted that this extension problem is not
at all easy; in fact, the main difficulty seems to be the extension
from $X$ to the boundary faces $F_i$. For this, it appears to be
almost necessary that $X$ and $M$ be products, and that the metrics
and boundary structures are globally rather similar to the ones
considered here.  Thus, in all that follows, we shall assume for
simplicity that $X$, $M$ and the distributions $\cD_i$ and $\cF_i$ are
products.

\subsection{Asymptotically product hyperbolic metrics}
\label{ssec:asympt-prod-hyperb}
It is always possible to construct a complete metric on the interior
of $M$ which is `weakly' product hyperbolic and with any given $(G_1×
G_2)$-conformal structure on $X$ as its prescribed conformal infinity.
In fact, we can just write down a formula which directly generalizes
(\ref{eq:thmet}): let $x_i$ be defining functions for the boundary
hypersurfaces $X_i \subset M_i$, and choose $\Im\, \bK_i$-valued $1$-forms
$\eta_i$ defining $\cD_i$ and compatible metrics $\gamma_i$ representing the
given conformal classes. In a neighbourhood of $X$ of the form
$(0,\epsilon)_{x_1} × (0,\epsilon)_{x_2} × X$, set
\begin{equation}
g_{\gamma_1,\eta_1,\gamma_2,\eta_2} = \lambda_1\big(\frac{dx_1^2}{x_1^2}+\frac{\gamma_1}{4x_1^2}+\frac{\eta_1^2}{4x_1^4}\big)
 + \lambda_2\big(\frac{dx_2^2}{x_2^2}+\frac{\gamma_2}{4x_2^2}+\frac{\eta_2^2}{4x_2^4}\big).
\label{eq:basicprod0}
\end{equation}
Slightly more generally, consider
\begin{equation}
g = g_{\gamma_1,\eta_1,\gamma_2,\eta_2} + k, \qquad k\in (x_1x_2)^\nu C^{2,\alpha}
\label{eq:basicprod}
\end{equation}
for some $\nu > 0$. (The norms and covariant derivatives are with
respect to the metric (\ref{eq:basicprod0}).)  We call any such $g$ a
weakly asymptotically product hyperbolic metric, and say that
$(\calD_1,\calD_2,[\gamma_1], [\gamma_2])$ is its conformal infinity. Just as
in the rank one setting, $g$ determines this conformal infinity.
Conversely, replacing $\gamma_i$ and $\eta_i$ in (\ref{eq:basicprod}) by any
other conformal representatives
\[
\widetilde{\gamma}_i = f_i \gamma_i, \qquad \widetilde{\eta}_i = f_i \eta_i,
\]
where the $f_i$ are strictly positive smooth functions on $X$, yields
a new metric
\[
\widetilde{g}  = 
\lambda_1\big(\frac{dx_1^2}{x_1^2}+\frac{f_1\gamma_1}{4x_1^2}+\frac{f_1^2\eta_1^2}{4x_1^4}\big)
 + \lambda_2\big(\frac{dx_2^2}{x_2^2}+\frac{f_2\gamma_2}{4x_2^2}+\frac{f_2^2\eta_2^2}{4x_2^4}\big).
\]
We claim that up to a diffeomorphism $\Phi$, $\widetilde{g}$ is
asymptotically equivalent to $g$. Indeed, if $\Phi^* \widetilde{x}_i =
\frac{x_i}{\sqrt{f_i}}$ and $\left. \Phi \right|_X = \mathrm{id}$, then
\[
\frac{d\widetilde{x}_i}{\widetilde{x}_i} = \frac{dx_i}{x_i}-\frac{df_i}{2f_i}, \qquad
\text{and}\qquad
\left|\frac{d\widetilde{x}_i}{\widetilde{x}_i}-\frac{dx_i}{x_i}\right|_g =  \calO(x_1+x_2),
\]
and hence 
\[
|\Phi^* \widetilde{g} - g|_g = \calO(x_1+x_2).
\]

A metric $g$ as in (\ref{eq:basicprod0}) is sometimes also called
weakly asymptotically Einstein, because of the
\begin{lemm} The metric (\ref{eq:basicprod0}) satisfies the estimate
\[
\left|\Ric^g + g\right|_g = \calO(x_1+x_2).
\]
\label{le:wae}
\end{lemm}
The proof is deferred to the next subsection, where the formalism for
the necessary calculations is developed.

As indicated earlier, we also define a narrower class of {\it strongly} 
asymptotically product hyperbolic (or strongly 
asymptotically Einstein) metrics, for which the Einstein tensor decays
uniformly near the entire boundary of $M$. The goal in the next few
subsections is to find the equations which the limiting values of
$g_1$ and $\gamma_2, \eta_2$ must satisfy at $x_2 = 0$ and the limiting
values of $g_2$ and $\gamma_1, \eta_1$ must satisfy at $x_1 = 0$, in order
that $g$ lie in this smaller class.

Later in the paper, in §~\ref{ssec:param-near-prod}, we shall also
define a class of `near product hyperbolic' metrics; these will be
defined by slightly different conditions, but we show there that any
strongly asymptotically product hyperbolic metric is of near product
type.

This profusion of similar names is indicative of the fact that
for metrics which are modelled by symmetric spaces of rank greater
than one, it is by no means clear what the precise conditions are
under which a metric should really be considered `asymptotically
symmetric'; each of the classes of metrics above has some claim to
this moniker in the product hyperbolic setting.

\subsection{Asymptotic curvature calculations}\label{ssec:asympt-curv-calc}
We now calculate the asymptotics of the Ricci curvature for a weakly asymptotically product 
hyperbolic metric (\ref{eq:basicprod}) on $M_1 × M_2$. In the course of this the proof 
of Lemma~\ref{le:wae} will emerge, as well as motivation for the extra conditions 
imposed on $g$ to warrant the name strongly asymptotically Einstein.

As before, assume that $M$, $X$ and the distributions $\cD_i \subset \cF_i$
are all of product type. The main calculations are local near the
boundary faces; to be definite we work in the region where $x_2 \to 0$,
and write $g$ in the form
\begin{equation}
g = g_1 + \lambda_2\left(\frac{dx_2^2}{x_2^2} + 
\frac{\gamma_2}{4x_2^2} + \frac{\eta_2^2}{4x_2^4}\right),\label{eq:1}
\end{equation}
where the two terms are metrics along horizontal and vertical slices,
$M_1 × \{p_2\}$ and $\{p_1\} × M_2$, respectively. We assume that
\begin{itemize}
\item the $x_2$ dependence is only what is written explicitly; in
  other words, $g_1$, $\eta_2$ and $\gamma_2$ are defined and smooth on $M_1
  × X_2$ and are independent of $x_2$;
\item $\eta_2$ is the pullback of a contact form from $X_2$, hence is
  independent of $M_1$;
\item $\gamma_2$ is a family of metrics on $\cD_2$ compatible with $\eta_2$ (and hence gives a 
$G_2$-conformal structure on each slice $\{p_1\}×X_2$).
\end{itemize}

The precise form of $g_1$ is not so important for the moment, but in
order to maintain consistency with (\ref{eq:basicprod}), we also
impose that in analogous coordinates near the boundary of $M_1$, $g_1
\sim \lambda_1 (dx_1^2/x_1^2 + \gamma_1/4x_1^2+ \eta_1^2/4x_1^4)$ as $x_1 \to 0$,
and that $\gamma_2$ converges to a representative of the specified
conformal class $[\gamma_2]$ as $x_1 \to 0$. However, these last conditions
do not enter into the immediate considerations.

Well-known formulæ due to O'Neill, cf.\ Proposition 9.36 in
\cite{Bes87}, express the Ricci curvature of a Riemannian submersion
in terms of the Ricci curvatures of the base and fibres and two
additional tensors: the second fundamental form $T$ of the fibres and
another tensor which measures the deviation of the horizontal
subspaces from being integrable. To adapt this to our setting, we
regard $M$ as a fibration $M_1 × M_2 \to M_1$. The two factors are
orthogonal, and the horizontal subspaces are integrable (with leaves
the $M_1$ slices, i.e.\ the submanifolds $M_1× \{q_2\}$), but this is
still not quite a Riemannian submersion because $g_1$ depends also on
$M_2$. In the curvature computations below, however, it behaves
asymptotically as $x_2 \to 0$ like a Riemannian submersion: the
negative powers of $x_2$ in all terms in $g_2$ add an extra $x_2$
factor to all derivatives in the $M_2$ directions.

We continue by defining the various quantities which appear in the
O'Neill formulæ, and developing some of their properties. The first is
the {\it second fundamental form} for $\gamma_2$.  This is the section $T$
of $T^*M_1 \otimes \Sym^2(\cD_2^*)$ defined by
\begin{equation}
\langle T(\xi_2,\zeta_2),\xi_1\rangle = - \frac{1}{2} (\cL_{\xi_1}\gamma_2)(\xi_2,\zeta_2)
\Leftrightarrow T = -\frac{1}{2}d^{M_1}\gamma_2.
\label{eq:def-T}
\end{equation}
(Here and later, a subscript $1$ or $2$ of a vector indicates the
factor to which it is tangent.)  For either of these expressions we
regard $\gamma_2$ as a section of $\Sym^2(\cD_2^*)$, which in turn
is a trivial bundle over each $M_1$ slice. For each $\xi_2\in \cD_2$ we
also set
\begin{equation}
T_{\xi_2} \in \mathrm{End}(\cD_2,TM_1), \qquad T_{\xi_2}^* \in \mathrm{End}(TM_1, \cD_2),
\label{eq:defTs}
\end{equation}
where the metrics $g_1$ and $\gamma_2$ are used to dualize.

There is still a freedom in the choice of the representative $\gamma_2 \in
[\gamma_2]$, but we now fix the normalization that the volume form
$dV^{\gamma_2}$ is constant in the $M_1$ directions. Consequently, the
mean curvature vector vanishes:
\[
\Tr^{\gamma_2} \, T=0.
\]

The trivial connection $d^{M_1}$ on $\cD_2$ is not compatible with the
metric, but to find one which is it suffices to add the second map in
(\ref{eq:defTs}); thus
\begin{equation}
\nabla = d^{M_1} + T, \qquad \left(\text{i.e.}\ \nabla_{\xi_1} \xi_2 = d^{M_1}_{\xi_1}\xi_2 + T^*_{\xi_2}\xi_1\right)
\label{eq:conn-slice}
\end{equation}
defines a unitary connection on $\cD_2$ over $M_1$. The divergence of
$T$ with respect to this connection is the bilinear form on $\cD_2$,
\[
(\delta^{M_1}T)(\xi_2,\zeta_2) = - \sum \langle(\nabla_{e_\alpha}T)(\xi_2,\zeta_2),e_\alpha\rangle,
\]
where $\{e_\alpha\}$ is an orthonormal frame for $TM_1$. The final
ingredient we need is the bilinear form $Q$ on $TM_1$ defined by
contracting the product of $T$ with itself in the
$\Sym^2(\cD_2^*)$ component with respect to $\gamma_2$,
\begin{equation}
Q(\xi_1 , \eta_1) := \langle T_{\xi_1},T_{\eta_1} \rangle_{\gamma_2}.
\label{eq:def-Q}
\end{equation}

Before proceeding, we derive the crucial first-order properties of
$T$.
\begin{lemm} Let $d_\nabla^{M_1}$ denote the exterior derivative on $M_1$
  coupled to the connection $\nabla$ on $\cD_2$. Then
\[
d_\nabla^{M_1}T=0.
\]
\end{lemm}
\begin{proof} By (\ref{eq:conn-slice}), $d_{\nabla}^{M_1}T = d^{M_1}T + T
  \land T$, where in the last term we regard the two factors as elements
of $T^*M_1 \otimes \mathrm{End}(\cD_2)$ and
$T^*M_1\otimes\Sym^2(\cD_2)$, respectively. The second expression
for $T$ in (\ref{eq:def-T}) gives $d^{M_1}T=0$, hence it suffices to
prove
\[
T\land T=0.
\]
This identity in turn is a direct consequence of the symmetry of the action 
\[
(u \cdot q)(x,y)=q(ux,y)+q(x,uy)=\gamma_2((uv+vu)x,y),
\]
of a symmetric endomorphism $u$ on a quadratic form
$q(x,y)=\gamma_2(vx,y)$ in the pair $(u,v)$.
\end{proof}
\begin{lemm}\label{lemm:bianchi}
On each slice $M_1 × \{z_2\}$ there is a Bianchi identity of the form
\[
\big(\delta^{M_1}Q+\frac{1}{2}d \, \Tr \, Q\big)_\xi = \langle\delta^{M_1}T,T^*\xi\rangle, \quad
\xi\in TM_1. 
\]

\end{lemm}
\begin{proof}
  Choose an orthonormal frame $\{e_\alpha\}$ for $TM_1$ and extend $\xi$ to
  a vector field on $M_1$ which is parallel with respect to $\nabla$ at
  some point $z_1$. Then, calculating at $z_1$,
\begin{align*}
(\delta^{M_1}Q)(\xi) &= - \sum_\alpha \nabla_{e_\alpha}Q(e_\alpha,\xi) \\
&= - \sum_\alpha \langle(\nabla_{e_\alpha}T^*)e_\alpha,T^*\xi\rangle + 
\langle T^*e_\alpha,(\nabla_{e_\alpha}T^*)\xi\rangle \\
&= \langle\delta^{M_1}T,T^*\xi\rangle - \sum_\alpha \langle T^*e_\alpha,(\nabla_\xi T^*)(e_\alpha)\rangle \\
&= \langle\delta^{M_1}T,T^*\xi\rangle - \frac{1}{2}d\Tr Q (\xi).
\end{align*}
The second equality uses that $d_\nabla^{M_1}T=0$. 
\end{proof}

There is also a second fundamental form $\bI_i$ and corresponding mean
curvature vector $N_i =\Tr^{g_i}\bI_i$ for each $M_i$ slice, $i=1,2$.
Note that $N_2$ is different from $\Tr^{\gamma_2}T$ (which we are assuming
is equal to $0$), since in the latter one only takes the trace in the
$\cD_2$ directions.

We can now state an exact formula for the Ricci curvature.
\begin{lemm} Let $(M=M_1 × M_2, g= g_1 + g_2)$ be a metric on $M$
  keeping the factors $M_1$ and $M_2$ orthogonal.  Let $\bI_i$ be the
  second fundamental form of $M_i$ and $N_i$ the mean curvature vectors.
  Then the Ricci tensor of $g$ is given by:
\begin{itemize}
\item 
\begin{equation}\begin{split}
\Ric^g(\xi_1,\zeta_1) = & \Ric^{g_1}(\xi_1,\zeta_1)- (\delta^{M_2}\, \bI_1)(\xi_1,\zeta_1) - 
\langle\bI_1(\xi_1,\zeta_1),N_1\rangle \\
& + (\delta^{M_1})^*\, N_2(\xi_1,\zeta_1) - \langle\bI_2^*\, \xi_1,\bI_2^*\, \zeta_1\rangle ,
\end{split}
\label{eq:ric-factor}
\end{equation}
with an analogous expression for the restriction of $\Ric^g$ to $TM_2$. 
Here $\delta^*$ is the symmetrization of the covariant derivative. 
\item 
\begin{equation}
\begin{split}
\Ric^g(\xi_1,\xi_2) = 
&\langle\delta^{M_1}\bI_1(\xi_1),\xi_2\rangle + \langle\nabla_{\xi_1}N_1,\xi_2\rangle \\
+ &\langle\delta^{M_2}\bI_2(\xi_2),\xi_1\rangle + \langle\nabla_{\xi_2}N_2,\xi_1\rangle.
\end{split}
\label{eq:ric-mixte}
\end{equation}
Here $\delta^{M_i}\bI_i$ is the divergence of $\bI_i$ regarded as a symmetric
2-tensor along $M_i$.
\end{itemize}
\end{lemm}
The derivations of these two formulæ are left to the reader.

As a first application, we have the

\begin{proof}[Proof of Lemma \ref{le:wae}] We will be applying (\ref{eq:ric-factor})
and (\ref{eq:ric-mixte}) with $M_i = (0,1)_{x_i} × X_i$. The second fundamental
form of the slices $M_1 × \{p_2\}$ is 
\[
\langle\bI_1,\xi_2\rangle  =-\frac{1}{2} \cL_{\xi_2}g_1 
=-\frac{\lambda_1}{2}(\frac{1}{4x_1^2}\cL_{\xi_2}\gamma_1+\frac{1}{4x_1^4}\cL_{\xi_2}\eta_1^2).
\]
If $\xi_2$ is a unit vector in $TX_2$, then 
\[
\frac{1}{4x_1^2}\cL_{\xi_2}\gamma_1+\frac{1}{4x_1^4}\cL_{\xi_2}\eta_1^2 = \calO(x_2),
\]
since $\cD_1$ does not depend on $M_2$. Hence on each slice
$M_1×\{x_2\}$, $\bI_1 = \calO(x_2)$, and the same is true for all its
derivatives. On these same slices one also has
\[
\Ric^{g_1} = -\lambda_1 g_1 + \calO(x_1) .
\]

Analogously, on the slices $\{p_1\} × M_2$, we have
\[
\bI_2 = \calO(x_1), \qquad \text{and} \qquad \Ric^{g_2} = -\lambda_2 g_2 + \calO(x_2) . 
\]
Inserting these in (\ref{eq:ric-factor}) and (\ref{eq:ric-mixte}) gives 
$\Ric^g= -g + \calO(x_1+x_2)$, as desired.
\end{proof}
\medskip

The main result of this subsection is the
\begin{lemm} Let $g$ be defined by (\ref{eq:1}), and suppose that 
\begin{equation}\label{eq:con-sli}\begin{cases}
& \delta^{M_1}T = 0 \\
& \Ric^{g_1}(\xi_1,\zeta_1)+\lambda_1\langle \xi_1,\zeta_1\rangle_{g_1} = 
\langle T^*\xi_1,T^*\zeta_1\rangle_{\gamma_2}
\end{cases}
\end{equation}
for all vectors $\xi_1,\zeta_1\in TM_1$. Then 
\[
\Ric^g=-g+ \calO(x_2).
\]
In particular, if $g$ is as in (\ref{eq:basicprod}) and the equations
(\ref{eq:con-sli}) are satisfied at both boundary hypersurfaces, then
$g$ is strongly asymptotically Einstein.
\end{lemm}

\begin{rema}
  Note that the right hand side of the Bianchi identity in Lemma
  \ref{lemm:bianchi} vanishes when the first equation in
  (\ref{eq:con-sli}) is satisfied.
\end{rema}

\begin{proof} We apply (\ref{eq:ric-factor}) and (\ref{eq:ric-mixte})
  as follows.  First, just as before, $\bI_1 = \calO(x_2)$, and
  hence $N_1=\calO(x_2)$ too.  On the other hand, the
  normalization on the volume form implies that $N_2=0$, thus
\begin{align*}
\Ric^g(\xi_1,\zeta_1) &= \Ric^{g_1}(\xi_1,\zeta_1) - \langle\bI_2^*\xi_1,\bI_2^*\zeta_1\rangle + 
\calO(x_2) , \\ 
\Ric^g(\xi_1,\xi_2) &= \calO(x_2) , \\
\Ric^g(\xi_2,\zeta_2) &= \Ric^{g_2}(\xi_2,\zeta_2) - (\delta^{M_1}\bI_2)(\xi_2,\zeta_2) +
\calO(x_2) .
\end{align*} 
{}From the formula (\ref{eq:1}) and the fact that $\eta_2$ is constant along $M_1$ slices, we get
$$
\bI_2=\frac{T}{x_2^2},
$$
and since the connection (\ref{eq:conn-slice}) on $\cD_2$ along $M_1$
is exactly the one induced by the Levi-Civita connection of $g$, the
result follows.
\end{proof}
To conclude the section, observe that the formulæ in the lemma
correspond exactly to the standard formulæ obtained for a Riemannian
submersion with integrable horizontal distribution \cite[proposition
9.36]{Bes87}, as expected from our claim that the asymptotic behaviour
when $x_2\to0$ is that of a Riemannian submersion.

\section{Extending the approximate solution to the codimension one boundary faces}
\label{sec:sef}

Let $(M=M_1 × M_2,g=g_1^0 + g_2^0)$ be a product of A$\bK$H Einstein
metrics, with conformal infinity $\frakc^0 =
(\cD^0_1,\cD^0_2,[\gamma_1^0],[\gamma_2^0])$ on $X = X_1 × X_2$.  As in the
last section, we consider deformations $\frakc$ of $\frakc^0$ (assumed
to be globally integrable in case either $X_1$ or $X_2$ is not simply
connected). According to Lemma \ref{lemm:small-def}, the pair of
distributions $\cD_1,\cD_2$ remains of product type on $X$, and we
then extend these by pullback to a pair of transverse distributions of
product type on all of $M$. Choose metrics $\gamma_i$ representing each of
the conformal classes $[\gamma_i]$. In the real or complex case, these may
depend on both factors of $X$, but we maintain the normalization so
that, still just over $X$, $dV^{\gamma_1}$ is independent of $X_2$, and
similarly $dV^{\gamma_2}$ is independent of $X_1$.  Based on the
calculations of §~\ref{ssec:asympt-curv-calc}, we now address the
problem of how to extend $(\gamma_1,\gamma_2)$ over the faces $M_1 × X_2$ and
$X_1 × M_2$ to obtain a metric which is strongly asymptotically
Einstein.

\subsection{Extension along boundary faces}
\label{sec:extens-along-bound}
We focus on the extension of $\gamma_1$ over $M_1 × X_2$, since
the other case is treated exactly the same.

Let us restate the problem more carefully. On the face $M_1 × X_2$, we
seek metrics $g_1$ on $TM_1$ and $\tilde{\gamma}_2$ on $\cD_2$
which solve the system (\ref{eq:con-sli}).  The solutions are
constrained by the requirements that $dV^{\tilde{\gamma}_2}$ is
independent of $M_1$, and that $(g_1,\tilde{\gamma}_2)$ is
asymptotic at $X = \del M_1 × X_2$ to the given $(G_1× G_2)$-conformal
structure in the sense that
\[
g_1 \sim \lambda_1\big(\frac{dx_1^2}{x_1^2}+\frac{\gamma_1}{4x_1^2}+\frac{\eta_1^2}{4x_1^4}\big),
 \qquad \tilde{\gamma}_2 \sim \gamma_2,
\]
as $x_1 \to 0$ (with an error term $\calO(x_1^\nu)$ for some $\nu > 0$).

\begin{lemm} Suppose that $\KK_2 = \HH$ or $\OO$. Then
  (\ref{eq:con-sli}) reduces to a single uncoupled equation on $M_1$
  which is simply the usual Einstein equation.
\end{lemm}
\begin{proof}
  Under this hypothesis, $\gamma_2$ is independent of the $X_1$ factor.
  Hence $T \equiv 0$ and the first equation in (\ref{eq:con-sli}) is
  satisfied. The second equation reduces to the uncoupled Einstein
  equation on $M_1$. By Proposition \ref{prop:deform-theory-abkh} we
  can extend the conformal class $[\gamma_1]$ on $X_1$ to an A$\bK_1$H
  Einstein metric $g_1$ on $M_1$; note that this is actually
  done parametrically, depending on $q_2 \in X_2$.
\end{proof}

When $M_2$ is real or complex, (\ref{eq:con-sli}) cannot be reduced in
this way, but fortunately, solutions can still be obtained near to the
standard one by perturbation methods.

We can now state and prove the main result of this section.
\begin{theo}\label{theo:equat-faces}
  Suppose $M=M_1 × M_2$ is a product of A$\bK$H Einstein spaces such
  that the $L^2$ nullspace for the linearized gauged Einstein operator
  on $M_1$ vanishes. Then, for any small globally integrable
  perturbation of the product $(G_1 × G_2)$-conformal structure on
  $X$, the system (\ref{eq:con-sli}) has a global solution
  $(g_1,\tilde \gamma_2)$ on the face $M_1 × X_2$ with the prescribed
  asymptotic behaviour at $X$, more precisely on each slice $M_1×\{q_2\}$,
\[
g_1 - \lambda_1 \big(\frac{dx_1^2}{x_1^2}+\frac{\gamma_1}{4x_1^2}+\frac{\eta_1^2}{4x_1^4}\big)
\in (x_1 x_2)^\nu \Lambda^{2,\alpha}, \qquad \tilde \gamma_2-\gamma_2 \in \calC^{2,\alpha}, 
\]
with smooth dependence with respect to $q_2$.
\end{theo}
\begin{proof}
  This proof is similar to that for Proposition
  \ref{prop:deform-theory-abkh}.  Consider the slice $M_1 × \{q_2\}$,
  and begin with the conformal structures $\gamma_1, \gamma_2$ on $X_1 × X_2$,
  with $dV^{\gamma_2}$ independent of $X_1$. Fix a smooth extension map
  assigning to $\gamma_1$ a metric $g_1$ on $M_1$ with
\[
g_1 = \lambda_1
\big(\frac{dx_1^2}{x_1^2}+\frac{\gamma_1}{4x_1^2}+\frac{\eta_1^2}{4x_1^4}\big)
+ k , \quad k\in (x_1 x_2)^\nu \Lambda^{2,\alpha}
\]
as $x_1 \to 0$. Here $\Lambda^{2,\alpha}$ is the geometric Hölder space on
$M_1 × \{q_2\}$. The weight $\nu$ is positive; we can fix $\nu=1$, but any
smaller value is possible. Recall that we already have
\[
\Ric^{g_1}+ g_1=\calO(x_1).
\]
As before, extend $\gamma_2$ by pullback on $M_1$, so that the
corresponding second fundamental form $T=-\frac12 d^{M_1}\tilde \gamma_2$
satisfies also $$T=\cO(x_1).$$

When $\KK = \RR$, we consider perturbations $\phi$ of $\gamma_2$ which fix
$dV^{\gamma_2}$, so the tangent space consists of trace-free symmetric
two-tensors, i.e.\ sections of
\[
\calS=\Sym^2_0(TX_2).
\]
When $\KK = \CC$, $\tilde \gamma_2(\cdot,\cdot)=d\eta_2(\cdot,J\cdot)$, where $J$ is an almost
complex structure; the perturbation $\phi$ must also remain compatible
with $d\eta_2$ on $\cD_2$, or equivalently is a deformation of $J$, so
the tangent space consists of trace-free $J$-skew-Hermitian symmetric
two-tensors, i.e.\ sections of
\[
\calS=\Sym^-_0(\cD_2).
\]

Consider elements
\[
h \in \Lambda^{2,\alpha}_\nu(M_1× \{q_2\}, \Sym^2(TM_1)), \qquad
\phi \in \Lambda^{2,\alpha}_\nu(M_1× \{q_2\}, \calS), 
\]
and assume that both $h$ and $\phi$ have sufficiently small norm.
Suppose that $(g_1+h,\tilde \gamma_2 + \phi)$ is a solution of the system
(\ref{eq:con-sli}) along $M_1$.

Denote by $T^\phi$ and $Q^\phi$ the second fundamental form and
corresponding quadratic form defined by $\gamma_2 + \phi$.  To break the
diffeomorphism invariance of the equation, we add the Bianchi gauge
condition
\[
B^{g_1}(h)=\big(\delta^{g_1}+\frac{1}{2}d\Tr^{g_1}\big)h=0.
\]
Thus we consider the system
\[
\Phi^{g_1,\tilde \gamma_2}(h,\phi)=\big(\Ric^{g_1+h} - \lambda (g_1 + h)-
Q^\phi+(\delta^{g_1+h})^*B^{g_1}(h), \delta^{g_1+h,\phi}T^\phi \big) .
\]
By Lemma \ref{lemm:bianchi}, any solution of this equation must also satisfy
\[
B^{g_1+h}(\delta^{g_1+h})^*B^{g_1}(h)=0.
\]
By the same argument as in the uncoupled case, cf.\ Proposition
\ref{prop:einst-equat-bianchi}, we conclude that $B^{g_1}h=0$.
Hence a solution of $\Phi^{g_1,\tilde \gamma_2}(h,\phi)=0$ is also a
solution of the original system (\ref{eq:con-sli}).

As in the proof of proposition \ref{prop:deform-theory-abkh}, it suffices to check that
\begin{multline*}
\Phi^{g_1,\tilde \gamma_2}: \Lambda^{2,\alpha}_\nu(M_1× \{q_2\}, \Sym^2(TM_1) \oplus \calS) \\
\longrightarrow \Lambda^{0,\alpha}_\nu(M_1× \{q_2\}, \Sym^2(TM_1) \oplus \calS)
\end{multline*}
is a $\calC^1$ mapping of Banach spaces, for $(h,\phi)$ of sufficiently
small norm, and furthermore, that its linearization
\begin{multline*}
\left. D\Phi^{g_1,\tilde \gamma_2}\right|_{(0,0)}: \Lambda^{2,\alpha}_\nu(M_1× \{q_2\}, \Sym^2(TM_1) \oplus \calS) \\
\longrightarrow \Lambda^{0,\alpha}_\nu(M_1× \{q_2\}, \Sym^2(TM_1) \oplus \calS)
\end{multline*}
is an isomorphism at the product metric. Since the linearization of
the Bianchi-gauged Einstein equation is $\nabla^*\nabla -2\overset{\circ}{R}$,
and $T=0$ at the product metric, this linearization decouples as
\[
\left. D\Phi\right|_{(0,0)}(\dot{h},\dot{\phi})=
\big((\nabla^*\nabla -2\overset{\circ}{R})\dot{h},\nabla^*\nabla \dot{\phi}\big)
\]
for the metric $g_1^0$. By the hypothesis on the vanishing of the
$L^2$ nullspace for $(M_1,g_1)$, the first component is an
isomorphism.  The second component $\nabla^*\nabla$ is the rough Laplacian,
and this is an isomorphism for weights $\nu \in (0,n+d-1)$.

The last statement comes from the smooth dependence of the solution
constructed by the inverse function theorem with respect to the
parameter $q_2$.
\end{proof}

\begin{coro}\label{cor:extens-along-bound-1}
  With the same hypotheses as in Theorem \ref{theo:equat-faces}, let
  $(g_1,\tilde \gamma_2)$ be the solution of the system (\ref{eq:con-sli})
  on the face $M_1 × X_2$.  Then the metric
\[
h = g_1 + \lambda_2\left(\frac{dx_2^2}{x_2^2} + \frac{\tilde \gamma_2}{4x_2^2} + \frac{\eta_2^2}{4x_2^4}\right),
\]
defined in some neighbourhood of the face $M_1 × X_2$ where $x_2 \ll
1$ on $M_1 × M_2$ satisfies:
\begin{enumerate}
\item $\Ric^h + h = O(x_2)$ uniformly on the closure of this face, and
  more precisely $ \Ric^h + h \in x_2 \Lambda^\alpha $;
\item when $x_1\to0$, then $h-h_0\in x_1\Lambda^{2,\alpha}$, where $h_0$ is the
  model metric (\ref{eq:basicprod0}) given by the formula
$$ h_0=\lambda_1\left(\frac{dx_1^2}{x_1^2} + \frac{\gamma_1}{4x_1^2}
              + \frac{\eta_1^2}{4x_1^4}\right) + 
       \lambda_2\left(\frac{dx_2^2}{x_2^2} + \frac{\gamma_2}{4x_2^2}
                     + \frac{\eta_2^2}{4x_2^4}\right).
$$
\end{enumerate}\end{coro}
\begin{proof} This is a direct consequence of the formulæ
  (\ref{eq:ric-factor}) and (\ref{eq:ric-mixte}) for $\Ric^g$, since
  the solution $(g_1,\gamma_2)$ depends smoothly on the parameter $q_2\in
  X_2$.
\end{proof}

We conclude this section with a comment about regularity. Exactly as
in Remark~\ref{re:reg} at the very end of
§~\ref{sec:einst-deform-theory}, the solutions obtained in
Theorem~\ref{theo:equat-faces} are polyhomogeneous at the boundaries
of the codimension one faces provided the $(G_1× G_2)$-conformal
infinity data on $X$ is smooth. The proof is identical to the one for
the uncoupled A$\KK$H Einstein equations.

\subsection{Strongly asymptotically Einstein metrics}
\label{sec:strongly-asympt-eins}
We are now ready to define, given any small deformation $\frakc$ of
the given $G_1×G_2$-conformal structure $\frakc^0$ on $X$, a global,
approximately Einstein metric $g$ on $M$. By Theorem
\ref{theo:equat-faces} and Corollary \ref{cor:extens-along-bound-1},
we extend the data $(\gamma_1,\gamma_2)$ on the two faces $F_1=M_1×X_2$ and
$F_2=X_1×M_2$ to get pairs $(g_1,\tilde\gamma_2)$ and $(\tilde\gamma_1,g_2)$
solving the system (\ref{eq:con-sli}) on each face. Therefore, the two
metrics on $M$,
$$ 
h_1=g_1+ \lambda_2\left(\frac{dx_2^2}{x_2^2} + \frac{~\tilde\gamma_2}{4x_2^2}
                     + \frac{\eta_2^2}{4x_2^4}\right), \qquad
h_2= \lambda_1\left(\frac{dx_1^2}{x_1^2} + \frac{\tilde\gamma_1}{4x_1^2}
              + \frac{\eta_1^2}{4x_1^4}\right) + g_2
$$
defined in the neighbourhoods $x_2 \ll 1$ and $x_1 \ll 1$ of
$F_1$ and $F_2$, respectively, satisfy
$$
h_1 - h_0 \in x_1 \Lambda^{2,\alpha}, \qquad h_2 - h_0 \in x_2 \Lambda^{2,\alpha},
$$
where
$$ 
h_0=\lambda_1\left(\frac{dx_1^2}{x_1^2} + \frac{\gamma_1}{4x_1^2}
              + \frac{\eta_1^2}{4x_1^4}\right) + 
       \lambda_2\left(\frac{dx_2^2}{x_2^2} + \frac{\gamma_2}{4x_2^2}
                     + \frac{\eta_2^2}{4x_2^4}\right)
$$
is the initial model metric (\ref{eq:basicprod0}). It remains to glue
$h_1$ and $h_2$ in the region $\{x_1\ll 1, x_2\ll 1\}$ where they both
exist. Choose some cut-off function $\chi$ such that $\chi(x)=1$ for
$x<1/2$ and $\chi(x)=0$ for $x>2$, and consider the metric
$$ 
h = \chi\big(\frac{x_2}{x_1}\big) h_1 + \Big(1-\chi\big(\frac{x_2}{x_1}\big)\Big) h_2 , 
$$
now defined in a neighbourhood of the whole boundary $F_1\cup F_2$. In the region $1/2<x_2/x_1<2$, 
all the derivatives of $\chi(x_2/x_1)$ remain bounded for the metric 
$\frac{dx_1^2}{x_1^2}+\frac{dx_2^2}{x_2^2}$, hence in that region $h-h_0\in x_1 \Lambda^{2,\alpha}$ 
(or equivalently $x_2 \Lambda^{2,\alpha}$), and therefore $ \Ric^h + h = \cO(x_1)$, 
or more precisely $\Ric^h+h\in x_1 \Lambda^{2,\alpha}$. Globally, in a neighbourhood of $F_1\cup F_2$,
\begin{equation}
\Ric^h + h \in x_1^\nu x_2^\nu \Lambda^\alpha , \qquad \nu=\frac{1}{2} .\label{eq:4}
\end{equation}

We now generalize this model. Fix $0 < \nu \leq 1/2$.
\begin{defi}\label{def:strongly-asympt-eins-1}
A metric $g$ on $M$ is \emph{strongly asymptotically Einstein} if it differs from 
the metric $h$ defined above by a term in $x_1^\nu x_2^\nu \Lambda^{2,\alpha}$.
\end{defi}
In particular, the Ricci curvature of any such metric satisfies (\ref{eq:4}). 
Note too that any other reasonably method of patching $h_1$ and $h_2$ together 
near the corner yields a metric $h'$ which is strongly asymptotically Einstein in 
this same sense.

\section{Generalized Laplacians on near-product hyperbolic spaces}\label{sec:gener-lapl-near}
We now discuss the construction of a parametrix for $L^g$ when $g$ is
a strongly asymptotically Einstein perturbation of a product
hyperbolic metric. Our goal is to show that $L^g$ is invertible
between two weighted Hölder spaces.  We do this in the following
steps. First, we analyze the Schwartz kernel of the inverse of $L^g$
when $(M,g)$ is exactly product hyperbolic using a contour integral
representation; we go on to obtain conormal bounds for this Schwartz
kernel on the `product hyperbolic double space' $M^2_{\ph}$. This
serves as an ansatz for the parametrix of $L^g$ when the metric $g$ is
weakly asymptotically product hyperbolic. We introduce a stronger
condition on $g$ of being near product hyperbolic, and show that under
this hypothesis we can construct a parametrix with Schwartz kernel
conormal on $M^2_{\ph}$ which is an inverse of $L^g$ up to a compact
error term. The final step is to show that each of these operators are
bounded between weighted Hölder spaces, which implies that $L^g$ is
Fredholm on these spaces. The fact that when $g$ is a product we have
an exact inverse for $L^g$ on $L^2$ which is bounded between these
Hölder spaces shows that $L^g$ is invertible between these spaces
then. The parametrix construction varies continuously with $g$, so we
conclude that $L^g$ remains invertible when $g$ is near-product
hyperbolic and sufficiently close to a product metric.

\subsection{The inverse of $L$ when $g$ is a product}
\label{ssec:pla}
We begin with an examination of the structure of the inverse of $L^g$
on $L^2(M)$ when $(M,g)$ is a product of A$\KK$H spaces. This is
mostly a review of the analysis in \cite{MazVas02}. We first present a
contour integral representation for the $L^2$ inverse $G$ of $L^g$
involving the resolvent families of the operators $L^{g_j}$ on each
factor; from this we deduce estimates for the pointwise off-diagonal
behaviour of the Schwartz kernel of $G$ using the analogous estimates
for the Schwartz kernels of the two constituent resolvents. We do not
make an effort to obtain the most precise pointwise estimates on $G$
here, but see \cite{MazVas02} and \cite{Huang} for more on this.

\subsection*{A representation formula for the Green function}
Let $(M,g)$ be a product of two A$\KK$H Einstein spaces, with
linearized gauged Einstein operator $L^g = \nabla^* \nabla -
2\overset{\circ}{R}$, acting on sections of 
$$\Sym^2 T^*M=\Sym^2T^*M_1 \oplus \left(T^*M_1\otimes T^*M_2\right) \oplus \Sym^2T^*M_2.$$
The operator $L^g$ preserves the three summands and acts by
\[
L^{g_1}+(\nabla^{M_2})^*\nabla^{M_2},\quad
(\nabla^{M_1})^*\nabla^{M_1}+(\nabla^{M_2})^*\nabla^{M_2}, \quad 
(\nabla^{M_1})^*\nabla^{M_1} + L^{g_2},
\]
respectively. In each of these three cases, it has the form
\[
L^g = L_1 \otimes I_{B_2} + I_{B_1} \otimes L_2 ,
\]
where $L_i$ acts on a Banach space $B_i$ of sections of a bundle $E_i$
on $M_i$, and $L^g$ acts on the completed tensor product $B_1
\widehat{\otimes} B_2$ of sections of $E=E_1\otimes E_2$ on $M$.

For the moment, let $B_i = L^2(M_i, E_i; dV_{g_i})$. Denote by
$R_i(\mu)$ the resolvent family $(L_i - \mu)^{-1}$. This is a
holomorphic family of bounded operators on $B_i$ for $\mu$ in the
resolvent set $\CC \setminus \mathrm{spec}\,(L_i)$; according to
Theorem \ref{th:resfam1},
\[
\mathrm{spec}\,(L_i) = \{\lambda_{ij}\}_{j=1}^{N_i} \cup [\alpha_i,\infty)
\]
for some $\alpha_i > 0$, with $\lambda_{ij} \in \RR$, $\lambda_{ij} < \alpha_i$.
Nondegeneracy of $(M_i,g_i)$ for the Einstein problem is the
assumption that $0 \notin \{\lambda_{ij}\}$.

The resolvent family of the product operator, $R(\lambda) = (L- \lambda)^{-1}$,
can be expressed as a sort of convolution of the resolvents on the two
factors. More precisely, it is proved in \cite{MazVas02} that
\begin{equation}
R(\lambda) = -\frac{1}{2\pi i} \int_{\Gamma_\lambda} R_1(\mu)\, R_2(\lambda-\mu)\, d\mu,
\label{eq:prodintres}
\end{equation}
where $\Gamma_\lambda$ is a contour lying in the common region of holomorphy of
the two factors in the integrand with ends converging linearly to $± i
\infty$ and such that the spectrum of $L_1$ lies entirely on one side and
the spectrum of $L_2$ lies entirely on the other side.  For us it
suffices to take $\lambda = 0$. In the simplest situation, neither $L_1$
nor $L_2$ have any negative eigenvalues, and in this case we take
$\Gamma_0 = i \RR$.  The general case, where one or the other does have
such eigenvalues, requires a slight modification to this formula.

Setting $R(0) = L^{-1}$, and $G$ its Schwartz kernel, if $f \in
\calC^\infty_0(M)$, then $u = Gf$ is the unique $L^2$ solution to the
equation $Lu = f$.

\subsection*{Estimates on $G$}
Fix local coordinates $z_j = (x_j,y_j)$ near a boundary point of
$M_j$, and denote by $K_j(z_j,z_j',\mu)$ the Schwartz kernel of
$R_j(\mu)$. Also, replace $\mu$ by $i\mu$. Then (assuming neither
operator $L_j$ has negative eigenvalues), the Schwartz kernel of
$G$ equals
\begin{equation}
G(z_1,z_2,z_1',z_2') = \frac{1}{2\pi } \int_{-\infty}^{\infty} K_1(z_1,z_1',i\mu)\, K_2(z_2,z_2',-i\mu)\, d\mu.
\label{eq:prodintker}
\end{equation}
We now obtain pointwise estimates for $G$ using the bounds on
$K_j(z_j,z_j',i\mu)$ in Proposition \ref{pr:expbounds}. Our goal is to
prove that $G$ lifts to a conormal distribution on a certain blowup of
$M^2$ which is conormal at all boundaries, polyhomogeneous along the
lifted diagonal, and has a leading polyhomogeneous term at the front
faces with conormal remainder.

First observe that the integral (\ref{eq:prodintker}) converges in the Banach space of bounded operators 
on $L^2(M;dV_g)$. This follows from the elementary estimate
\[
||R_j(i\mu)||_{{\mathcal B}(L^2(M_j))} \leq \frac{1}{|\mu|},
\]
which is a direct consequence of the spectral theorem for the
selfadjoint operator $L_j$. This already uniquely specifies $G$ as an
element of ${\mathcal D}'(M^2)$.

Next, introduce the product-hyperbolic double space
\[
M^2_{\ph} = (M_1)^2_{\KK_1\Theta} × (M_2)^2_{\KK_2\Theta};
\]
this manifold with corners is simply the product of the $\KK_j \Theta$
double spaces of the two factors.  It has six codimension one boundary
faces: the two front faces, which are the boundary hypersurfaces
intersecting the lifted diagonal $\mathrm{diag}_{\ph}$:
\[
\mathrm{ff}_1  = \mathrm{ff}((M_1)^2_{\KK_1 \Theta}) × (M_2)^2_{\KK_2\Theta} \qquad \text{and} \qquad
\mathrm{ff}_2 = (M_1)^2_{\KK_1 \Theta} × \mathrm{ff}((M_2)^2_{\KK_2\Theta}),
\]
and the four side faces, which are the products of the side faces of
one factor with the interior of the other factor. We denote the side
faces by $B_{10,j}$ and $B_{01,j}$, $j = 1,2$, and for uniformity also
write $\mathrm{ff}_j = B_{11,j}$; the $j$ signifies that the face in
question comes from a boundary face in the $M_j$ factor. Defining
functions for any one of these faces will be written $\rho_{pq,j}$,
$j=1,2$ and $pq = 10, 01, 11$.

We shall prove that $G$ is conormal at all boundary faces of
$M^2_{\ph}$ and along the lifted diagonal. The conormal estimates here
are not quite sharp; there is a subtle cancellation in
(\ref{eq:prodintker}), explored more carefully in \cite{MazVas02},
which leads to vanishing rates which are slightly better (by a
logarithmic factor), but this is not needed here.

\begin{prop} The lift of the Schwartz kernel $G$ to the space
  $M^2_{\ph}$ is an element of
  $\calA^\sigma_{\ff}(M^2_{\ph},\mathrm{diag}_{\ph})$, where $\sigma =
  (\sigma_{pq,j})$ is the index set with $\sigma_{10,j}=\sigma_{01,j} = n_j +
  d_j-1$, $\sigma_{11,j} = 0$, and where $\ff = \ff_1\cup \ff_2$.
\end{prop}

\begin{proof}
  To separate out the contributions from the near-diagonal parts of
  each factor, we use standard results concerning the symbol calculus
  with spectral parameter, cf.\ \cite{Shu01}. Write $K_j = A_j + B_j$
  where $A_j$ contains the full diagonal singularity, is supported
  near the lifted diagonal in $(M_j)^2_{\KK_j \Theta}$, and is smooth
  across the front face of this double-space, while $B_j \in
  \calA^{\sigma_j}_\ff((M_j)^2_{\KK_j \Theta})$, $\sigma_j = ((n_j + d_j-1)/2 + 
\delta_0^{\KK_j},(n_j+d_j-1)/2 + \delta_0^{\KK_j})$.

  The integral of $A_1(i\mu)A_2(-i\mu)$ is a distribution supported in a
  neighbourhood of $\mathrm{diag}_{\ph}$ which intersects only the
  front faces of $M^2_{\ph}$ but not the side faces. It may be
  estimated directly (just as for the analogous computation on a
  product of two compact manifolds) using the oscillatory integral
  representations in the conormal bundle of $\mathrm{diag}_{\ph}$, and
  the result is clearly smooth up to the front faces. This is the
  Schwartz kernel of a pseudodifferential operator of order $-2$, as
  required.

  Next, the norm of $A_j$ as a map on any fixed weighted Sobolev or
  weighted Hölder space decays like $1/|\mu|$, cf. \cite{Shu01}. (Note
  that this is equivalent to the situation in the compact case because
  of the support properties of $A_j$.) From this we deduce that the
  integrals of $A_1(i\mu)B_2(-i\mu)$ and $A_2(i\mu) B_1(-i\mu)$ are smooth
  on the interior of $M^2_{\ph}$. By taking any number of derivatives
  with respect to $b$-vector fields on the $B$ factor, we obtain the
  correct conormal estimates too.

  As for the integral of $B_1(i\mu)B_2(-i\mu)$, divide the contour
  $\Gamma_0$ into a compact portion, where $|\mu| \leq \mu_0$, and its
  remaining noncompact ends $|\mu| \geq \mu_0$. The conormal estimates on
  $B_j$ hold locally uniformly in $\mu$, so in particular
\[
\int_{|\mu| \leq \mu_0} B_1(z_1,z_1',i\mu) B_2(z_2,z_2',-i\mu)\, d\mu \in \calA^{\sigma}_{B_{11}}(M^2_{\ph}).
\]
On the other hand, using the exponential bounds on $K_j$ from
(\ref{pr:expbounds}) we immediately deduce that the integral over
$|\mu| > \mu_0$ satisfies the same $\calC^0$ bound. The estimates for
higher tangential derivatives required to check conormality are
obtained in exactly the same way, using the corresponding pointwise
bounds for the higher tangential derivatives on each of the two
factors.

This completes the estimation of the Schwartz kernel $G$. 
\end{proof}

\subsection*{Modifications when either factor has negative point spectrum}
If $\inf \mathrm{spec}\,(L_1)$ or $\inf \mathrm{spec}\,(L_2)$
are negative, then (\ref{eq:prodintker}) needs to be altered slightly.
To understand this, begin by noting that (\ref{eq:prodintres}) remains
valid when $\lambda = i\epsilon$ and the contour separates the sets $i\epsilon -
\mathrm{spec}\,(L_1)$ and $\mathrm{spec}(L_2)$ in two
different half-planes. In fact, we take the contour to be the union of
the vertical rays $[2i\epsilon, i\infty)$, $(-i\infty, -i\epsilon]$ and two long thin
half-ellipses, one in the first quadrant with minor axis connecting
$\frac12 i \epsilon$ and $2i\epsilon$, and the other in the third quadrant, with
minor axis connecting $-\frac12 i \epsilon$ and $\frac12 i\epsilon$.  Now let $\epsilon
\searrow 0$; we arrive at the formula that $(L_1 + L_2)^{-1}$
is equal to the sum
\begin{multline*}
-\frac{1}{2\pi i }\int_{-i \infty}^{i\infty} R_1(-\mu)R_2(\mu)\, d\mu  \\ +
\sum_{\mu_j^{(1)}<0} P_1(-\mu_j^{(1)}) R_2(\mu_j^{(1)}) + 
\sum_{\mu_j^{(2)}<0} R_1(-\mu_j^{(2)}) P_2(\mu_j^{(2)}).
\end{multline*}
Here $\mu_j^{(i)}$ are the (negative) eigenvalues of $L_i$, and
$P_i(\mu_j^{(i)})$ is the orthogonal projection onto the corresponding
eigenspace in the $i^{\mathrm{th}}$ factor.

The main result in the last subsection, that the Schwartz kernel of
$(L^g)^{-1}$ is conormal on $M^2_{\ph}$, clearly remains valid. The
analysis of the main term, which is the first summand in this
expression, is exactly the same as before. The remaining terms are
much simpler to analyze since each is a simple tensor product of one
term which is polyhomogeneous on one factor and a finite rank
polyhomogeneous term on the other factor.

\subsection{The parametrix in the near-product case}\label{ssec:param-near-prod}
If $g = g_1 + g_2$ is a nondegenerate product hyperbolic metric, we
have established that the Schwartz kernel of the $L^2$ inverse of
$L^g$ has only conormal singularities on the product hyperbolic double
space. We now show that a similar structure theorem holds for a
parametrix for $(L^g)^{-1}$, i.e.\ an approximate inverse up to
compact error terms, when $g$ satisfies a condition we shall call
being `near product-hyperbolic'; the precise definition is given
later. This parametrix $H$ will have Schwartz kernel in
$\Psi^{-2}_{\ph}(M) + \calA^{\sigma}_{\ff}(M^2_{\ph})$. The construction of
$H$ uses the structure of $(L^g)^{-1}$ in the product case as an
ansatz and a tool. We describe this now, deferring the statement of
the final result to a theorem at the end of this subsection.

The procedure is very much the same as in the A$\KK$H case. The
Schwartz kernels of general $\ph$ pseudodifferential operators are
defined to be distributions on $M^2_{\ph}$ which are conormal at the
boundaries and with a polyhomogeneous expansion at the diagonal. There
is a small calculus of operators with Schwartz kernels supported near
the lifted diagonal, and a large calculus which also admits operators
with Schwartz kernels conormal up to all faces, and with positive
vanishing order at the front faces. The corresponding decomposition of
$H$ is written $H_1 + H_2$.

Let $g$ be a weakly asymptotically product hyperbolic metric on $M$.
The initial approximation $H_1$ to $(L^g)^{-1}$ is obtained via the
standard elliptic parametrix construction using the symbol calculus on
the conormal bundle of $\mathrm{diag}_{\ph}$ in $M^2_{\ph}$. This uses
the uniform invertibility of the $\ph$-symbol of $L^g$ up to the two
front faces, which is valid for any weakly asymptotically product
hyperbolic metric. If $g$ is only polyhomogeneous at the boundaries of
$M$, then $H_1 \in
\calA^{\tau}_{\mathrm{ff}}(M^2_{\ph},\mathrm{diag}_{\ph})$, where
$\tau_{10,j}=\tau_{01,j} = \infty$ and $\tau_{11,j}=0$, $j=1,2$.

Set $L^g H_1 = I - Q$. By construction, $Q \in
\calA^{\tau}_{\ff}(M^2_{\ph})$. The correction term $H_2$ is chosen so
that $L (H_1 + H_2) = I - Q_1$, or equivalently, $LH_2 = Q - Q_1$,
where $H_2 \in \calA^{\sigma}_{\ff}(M^2_{\ph}, \mathrm{diag}_{\ph})$ and
$Q_1 \in \calA^{\tilde{\sigma}}(M^2_{\ph})$; here $\sigma_{11,j} = 0$,
$\tilde{\sigma}_{11,j} > 0$, and all other $\sigma_{pq,j} = \tilde{\sigma}_{pq,j}
= n_j + d_j -1 $ for $pq \neq 11$. The vanishing of the Schwartz kernel
of $Q_1$ to some positive order at the front faces yields its
compactness on weighted Hölder spaces.

To determine $H_2$, restrict the equation $L^g LH_2 = Q-Q_1$ to each
of the two front faces. Set
\[
N_j(L^g) =  \left. L^g \right|_{\mathrm{ff}_j},\ N_j(Q_0) = \left. Q_0 \right|_{\mathrm{ff}_j}, \quad j = 1,2.
\]
Then we must solve the two equations $N_j(L) H_j = N_j(Q_0)$. The
right hand side is a smooth function on the interior of $\ff_j$; it
vanishes to infinite order at the intersection with the side faces,
and is smooth up to the other front face, or at least, the dependence
in this direction is exactly as regular as the metric $g$ near the
corner.

Let us begin by analyzing the structure of $N_j(L)$. Fixing $j=2$ to
be definite, the front face $\mathrm{ff}_2$ is a product
$(M_1)^2_{\KK_1\Theta} × \ff(M_2)$. The second factor, the front face of
$(M_2)^2_{\KK_2\Theta}$, is a fibration with base space the diagonal of
$(\del M_2)^2$ and each fibre naturally identified with the hyperbolic
space $\KK_2H^{m_2}$. The lift $L^g$ from the left factor of $M$ to
$M^2_{\ph}$ acts on the `left factor' of $M_j$ in each $\KK_j \Theta$
double space. Its restriction to $\ff_2$, which makes sense since it
acts tangentially to that boundary, is a sum of derivatives of two
types: some act on the left factor of $M_1$ in $(M_1)^2_{\KK_1\Theta}$ and
the others act on the $\KK_2 H^{m_2}$ fibres of $\ff(M_2)$. In
particular, the dependence on all other variables -- namely, in the
right factor of $M_1$ and in the diagonal of $(\del M_2)^2$ -- is
purely parametric. The key assumption is that $N_2(L^g) = L_{2,1} +
L_{2,2}$, where $L_{2,1}$ is an operator acting on the left factor of
$M_1$ (with coefficients depending parametrically on the right factor
of $M_1$ and the diagonal of $(\del M_2)^2$) and $L_{2,2}$ is the
gauged linearized Einstein operator or the rough Laplacian on
$\KK_2H^{m_2}$. Similarly, $N_1(L^g) = L_{1,1} + L_{1,2}$, where these
summands have analogous properties.

Let us investigate when this condition holds. Using interior
coordinates $z_1$ on $M_1$ and coordinates $z_2 = (x_2,y_2)$ near the
boundary on $M_2$, we have near $x_2 = 0$,
\[
L^g = \sum a^{(1)}_{\alpha}(z_1, y_2)\del_{z_1}^\alpha + \sum a^{(2)}_{j,\beta,\gamma}(z_1,y_2)(x_2 \del_{x_2})^j
(x_2 Y)^\beta (x_2^2 Z)^\gamma + \calO(x_2);
\]
where $(Y_1, \ldots, Y_{d_2(m_2-1)})$ and $(Z_1,\ldots, Z_{d_2-1})$ are bases
of sections of $\cD_2$ and its complement, respectively, in $TX_2$,
and we are using obvious multi-index notation. Notice that we are
using the product structure $M_1 × \del M_2$ for this boundary face.
Near product hyperbolicity is simply the condition that there are {\it
  no} cross-terms at $x_2 = 0$, and similarly at $x_1 = 0$. Note that
by naturality of generalized Laplacians, the second summand here, at
$x_2 = 0$, must be the linearized gauged Einstein operator on $\KK_2
H^{m_2}$.

\begin{lemm}
  Let $g$ be a strongly asymptotically product hyperbolic metric, as
  constructed in §~\ref{sec:strongly-asympt-eins}. Then $g$ is near
  product hyperbolic.
\end{lemm}
\begin{proof} By the global rigidity assumption for the conformal
  infinity data, we are fixing the product structure on each boundary
  hypersurface. Consider the construction of $g$ at $M_1 × \del M_2$.
  The preliminary extensions of the conformal infinity data $(\gamma_1,
  \eta_1)$ and $\gamma_2/x_2^2 + \eta_2/x_2^4$, regarded as metrics on the
  bundle $TM_1 \oplus {}^{\KK_2\Theta}T_{\del M_2}M_2$ over this face, make
  these subbundles orthogonal to one another. The correction terms $h$
  and $\phi$ are chosen to preserve this orthogonality.
\end{proof}

Now, suppose that $L_g$ has this near product form. Consider its lift
to $M^2_{\ph}$; near $\ff_2$ we replace $(z_2,\tilde{z}_2)$ with
projective coordinates, so $y_2 = \tilde{y}_2 + \calO(\tilde{x}_2)$
and the derivatives with respect to $x_2$, $Y_i$ and $Z_j$ in the
second summand all become tangential to $\KK_2 H^{m_2}$. Setting
$\tilde{x}_2=0$ gives $N_2(L)$. The coefficients $a^{(1)}_\alpha$ then
depend on $(z_1,\tilde{y}_2)$, but not the second hyperbolic space
factor; the coefficients $a^{(2)}_{j,\beta,\gamma}$ potentially depend on all
variables on this face, although the derivatives just act tangentially
to the hyperbolic spaces, but $L_{2,2}$ is simply the linearized
gauged Einstein operator or the rough Laplacian on $\KK_2H^{m_2}$, and
hence is independent of $z_1$ and $\tilde{y}_2$.  In other words, if
$g$ is any near product hyperbolic metric, in particular a strongly
asymptotically Einstein metric from Theorem~\ref{theo:equat-faces},
then the restriction of $N_j(L^g)$ to the fibres of the front face
$\ff_j$ is of product type with respect to the product structure on
each fibre. A corresponding statement is true for any generalized
Laplacian $P^g$ associated to $g$.

Having established this, we now proceed to solve the two normal
equations using the results of §~\ref{ssec:pla}. More specifically,
since $N_j(L^g)$ is of product type, the results of §~\ref{ssec:pla}
show that $N_j(H_2) = N_j(L^g)^{-1}N_j(Q_0)$ is conormal on $\ff_j$
and smooth up to the corner $\ff_1 \cap \ff_2$. To verify this last
statement, note that the product analysis of §~\ref{ssec:pla} gives
conormality on each slice $M_1 × \KK_2 H^{m_2}$ (when $j=2$, for
example), but the dependence on the other variables, i.e. the `right'
factor of $M_1$ in $(M_1)^2_{\KK_1\Theta}$ and the base space variable in
$\ff((M_2)^2_{\KK_2\Theta}$, are parametric, and the solution is as smooth
as the dependence on these parameters.  As the variable in the right
factor of $M_1$ converges to the boundary, the space on which
$N_2(L^g)$ acts converges locally uniformly to the product of the two
hyperbolic spaces; the uniform estimates on $N_2(Q_0)$ and on
$N_2(L^g)^{-1}$ show that the solutions converge smoothly too, and
that the limiting solution is conormal on the product of the radial
compactifications of the hyperbolic spaces.  In particular, the corner
$\ff_1 \cap \ff_2$ is a bundle over $\del M_1 × \del M_2$ with fibre the
product of hyperbolic spaces, and the limit of the solutions from
either side, $\ff_1$ or $\ff_2$ are the same.  Therefore
\[
\left. N_1(H_2) \right|_{\ff_1 \cap \ff_2} = \left. N_2(H_2) \right|_{\ff_1 \cap \ff_2}.
\]

This argument proves the existence of conormal functions on the two
front faces which are compatible at the corner. Hence we may extend
these to a Schwartz kernel $H_2$ in the interior of $M$ which is
conormal to all boundary faces, $H_2 \in
\calA^{\sigma}_{\ff}(M^2_{\ph})$, and which restricts to $N_j(H_2)$ at
$\ff_j$.

We have now constructed $H = H_1 + H_2$ which satisfies $L^g H = I -
Q_1$, $Q_1 \in \calA^{\tilde{\sigma}}(M^2_{\ph})$.  A left parametrix is
obtained by taking adjoints. This completes the parametrix
construction when $g$ is near product hyperbolic and polyhomogeneous.

As explained earlier, we can certainly restrict to studying
polyhomogeneous metrics with smooth conformal infinity data but we
comment briefly on how to extend this proof to the case where $g$ is
the sum of a polyhomogeneous A$\KK$H metric and a perturbation term $k
\in x_1^{\nu_1}x_2^{\nu_2}\Lambda^{2,\alpha}_{\KT}$, where $0 < \nu_j < n_j +
d_j-1$. (Note that this is the regularity for the strongly
asymptotically Einstein metrics.)

The symbol calculus step goes through immediately for symbols with
this regularity. Indeed, the perturbation term $k$ appears only in the
parametric dependence along the diagonal and does not occur in the
leading terms near the front faces, so $H_1$ decomposes as a principal
polyhomogeneous term and another, vanishing to some positive order at
these faces but which is only bounded in Hölder norm. The restriction
to the front faces of the error term $Q = I - LH_1$ does not depend on
these lower order terms, and the normal operators $N_j(L)$ are
independent of them too, except as above in their dependence on
parameters, which means that the second step carries through exactly
as before.

We have proved the 
\begin{prop} Let $g$ be a strongly asymptotically Einstein metric as
  constructed in §~\ref{sec:strongly-asympt-eins}. Then there is a
  parametrix $H = H_1 + H_2$ for $L^g$, with $H_1 \in
  \calA^{\sigma}_{\ff}(M^2_{\ph}, \mathrm{diag}_{\ph})$ and $H_2 \in
  \calA^\sigma_{\ff}(M^2_{\ph})$. The error term $Q_1 = I - L(H_1 + H_2)$
  is in $\calA^{\tilde{\sigma}}_{\ff}(M^2_{\ph})$.
\end{prop}

\subsection{Function spaces and mapping properties}
We now show that the pseudodifferential inverses or parametrices for
$L^g$, when $g$ is product, or near product, hyperbolic, are bounded
on weighted `geometric' Hölder spaces associated to the metric $g$.

The definition of these Hölder spaces is essentially identical to the
one in §~\ref{sec:rola}, the only difference when the weight parameters are $0$
being that to define $\Lambda^{k,\alpha}_{\ph}$ when $k > 0$, we allow
differentiations by arbitrary smooth vector fields which are locally
given as finite combinations of smooth products of elements of
$\calV_{\KK_1\Theta}(M_1) × \calV_{\KK_2\Theta}(M_2)$. If $\nu_1$ and $\nu_2$
are any two weight parameters, we also define
$x_1^{\nu_1}x_2^{\nu_2}\Lambda^{k,\alpha}_{\ph}(M)$ in the obvious way.

If $\KK_1 = \HH$, then the Hölder spaces corresponding to two
different near product hyperbolic metrics $g$ and $g'$ are not
quasi-isometric to one another, unless the corresponding distributions
$\cD_1$ and $\cD_1'$ on $\del M_1 × M_2$ are diffeomorphically
equivalent.

\begin{prop} Let $H$ be an operator which has its Schwartz kernel in
$\calA^{\sigma}_{\ff}(M^2_{\ph}, \mathrm{diag}_{\ph})$, where $\sigma$ is
the usual weight family with $\sigma_{11,j} = 0$ and all other
$\sigma_{pq,j} = (n_j + d_j-1)/2 + \delta_0^{\KK_j}$, 
and with pseudodifferential order $-\ell$ along $\mathrm{diag}_{\ph}$.  
Let $0 < \nu_j < (n_j + d_j-1)/2 + \delta_0^{\KK_j}$. Then
\[
H: x_1^{\nu_1}x_2^{\nu_2}\Lambda^{k,\alpha}_{\ph}(M) \longrightarrow
x_1^{\nu_1}x_2^{\nu_2}\Lambda^{k+\ell,\alpha}_{\ph}(M)
\]
defines a bounded map. If $Q$ has Schwartz kernel in
$\calA^{\tilde{\sigma}}(M^2_{\ph})$ where $\tilde{\sigma}_{11,j} > 0$ but all
other weights are the same as for $\sigma$, then $Q$ defines a compact
mapping between these same spaces.
\end{prop}
\begin{proof} As before, write $H = H_1 + H_2$. The boundedness of
  $H_1$ follows from the standard boundedness of pseudodifferential
  operators of on Hölder spaces as well as the fact that the support
  conditions on this Schwartz kernel means that this part of $H$ only
  spreads supports by a fixed amount.

  The proof of the $\calC^0$ bound for $H_2f$ when $f$ is in this
  weighted Hölder space proceeds by a direct and elementary estimation
  of this integral. Any higher derivative of $H_2 u$ with respect to a
  $\KK_1\Theta$ vector field on the first factor or a $\KK_2\Theta$ vector
  field on the second factor, or an iterated combination of such
  vector fields, is handled by simply noting that any of these vector
  fields applied to $H_2$ gives a Schwartz kernel of exactly the same
  form, with the same orders of vanishing at all the side faces.

  Finally, to prove that $Q$ is a compact operator, the extra
  vanishing at the front faces gives that
\[
Q: x_1^{\nu_1}x_2^{\nu_2}\Lambda^{s_1,\alpha}_{\ph}(M) \longrightarrow 
x_1^{\nu_1 + \e}x_2^{\nu_2 + \e }\Lambda^{s_2,\alpha}_{\ph}(M)
\]
for any $s_1, s_2 \in \RR$. This range space on the right includes
compactly into the domain space on the left (provided $s_2 > s_1$),
using the Arzela-Ascoli theorem.
\end{proof}

\begin{theo} Let $g$ be a near product hyperbolic metric which is
  sufficiently close to a nondegenerate product hyperbolic metric
  $g_0$ on $M_1 × M_2$. Then
\[
L^g: x_1^{\nu_1}x_2^{\nu_2}\Lambda^{2,\alpha}_{\ph}(M, S^2T^*M) \longrightarrow 
x_1^{\nu_1}x_2^{\nu_2}\Lambda^{0,\alpha}_{\ph}(M,S^2 T^*M)
\]
is an isomorphism.
\label{th:isohold}
\end{theo}
\begin{proof} If $g$ is any near product hyperbolic metric, then we
  have constructed a parametrix $H^g$ for $L^g$ so that $L^g H^g = H^g
  L^g = I - Q^g$.  We have also proved that $H^g$ is a bounded
  operator between the space on the right and the space on the left
  above, and that $Q^g$ is compact between these same spaces. This
  proves that the map $L^g$ is Fredholm.

  The construction of $H^g$ depends continuously on the metric (with
  respect to some sufficiently strong topology).  Since $Q^g = 0$ when
  $g = g_0$, we can make the norm of $Q^g$ as small as desired when
  $g$ is sufficiently close to $g_0$, which implies that $L^g$ is
  invertible.
\end{proof}

\section{Solving for the Einstein metric}
\label{sec:gedt}
The remainder of the proof of the Main Theorem proceeds very much as
in the analogous arguments in §~\ref{sec:einst-deform-theory} and
§~\ref{sec:sef}.

Let $(M=M_1× M_2,g = g_1 + g_2)$ be a product of A$\KK$H Einstein
metrics. Let $L^g$ and $L^{g_i}$ be the linearized gauged Einstein
operators for $g$ and the two component metrics $g_i$. We assume that
\[
0 \notin \mathrm{spec}\,(L^g) \cup \mathrm{spec}\,(L^{g_1}) \cup \mathrm{spec}(L^{g_2}),
\]
and that at least one of $\KK_1$ or $\KK_2$ is not equal to $\HH$ or
$\OO$.  Denote by $\frakc$ the $(G_1 × G_2)$-conformal infinity data
on $X = X_1 × X_2 = \del M_1 × \del M_2$.

\begin{theo}
  Under all these conditions, let $\frakc'$ be any globally integrable
  $(G_1 × G_2)$-conformal infinity data on $X$ which is $\calC^\infty$ and
  sufficiently close (in $\calC^{2,\alpha}$ norm) to $\frakc$.  Then there
  is a unique near product hyperbolic Einstein metric $g'$ which is
  close to the near product hyperbolic asymptotically Einstein metric
  $\tilde{g}'$ with conformal infinity data $\frakc'$ constructed in
  §~\ref{sec:strongly-asympt-eins}.
\end{theo}
\begin{proof} Let $\tilde{g}'$ be the strongly asymptotically Einstein
  metric constructed in §~\ref{sec:strongly-asympt-eins} with
  conformal infinity data $\frakc'$.  Write $g' = \tilde{g}' + k$,
  where
\[
k \in x_1^{\nu}x_2^{\nu}\Lambda^{2,\alpha}_{\ph}(M, \Sym^2(T^*M)).
\]
Write the gauged Einstein equation as
\[
N^{\tilde{g}'}(k) = \Ric(\tilde{g}' + k) + (\tilde{g}' + k) + 
(\delta^{\tilde{g}'+k})^*B^{\tilde{g}'}(k) = 0.
\]
The linearization at $k=0$ is the generalized Laplacian $L^{\tilde{g}'}$. 

We have that 
\[
N^{\tilde{g}'}: x_1^{\nu_1}x_2^{\nu_2}\Lambda^{2,\alpha}_{\ph}(M, \Sym^2(T^*M)) \longrightarrow
x_1^{\nu_1}x_2^{\nu_2}\Lambda^{0,\alpha}_{\ph}(M, \Sym^2(T^*M))
\]
is a $\calC^1$ map from a neighbourhood of $0$ in the domain space, and that 
$L^{\tilde{g}'}$ is a bounded linear map between these same two spaces.

According to Theorem~\ref{th:isohold}, if $\frakc'$ is sufficiently
close to $\frakc$, this linearization is an isomorphism. Furthermore,
the norm of its inverse is bounded away from zero, uniformly as
$\frakc' \to \frakc$.  The inverse function theorem implies that there
is a unique solution $k$ to $N^{\tilde{g}'}(k) = 0$ with $k$ near $0$.
\end{proof}

\appendix
\section{Osculating quaternionic coordinates}
\label{app:qcs}
In this brief appendix we prove Lemma \ref{lem:coord-q}. This uses an
idea close to that used to find normal coordinates in Riemannian
geometry, and should be a general fact for all the so called
``parabolic geometries''.

Fix a quaternionic contact structure on $Y^{4m-1}$ and any metric in
this conformal class on the distribution $\calD$. We use the
Tanaka-Webster type connection $\nabla$ from \cite[chapter 2]{Biq00}, see
also \cite{Duc04,Duc06} for the special case of dimension $7$. The contact
distribution $\cD$ is the kernel of three $1$-forms
$(\eta_1,\eta_2,\eta_3)$, and has a privileged supplementary subspace
generated by three ``Reeb vectors fields'' $R_1,R_2,R_3$ which are
uniquely specified by the conditions $\eta_i(R_j)=\delta_{ij}$ and
$(i_{R_i}d\eta_j+ i_{R_j}d\eta_i)|_\cD=0$ (in dimension $7$ a weaker
condition is placed on these). This connection satisfies:
\begin{itemize}
\item $\nabla$ preserves the distribution $\cD$ and the quaternionic
  structure on $\cD$;
\item the torsion $T^\nabla$ of two horizontal vectors $X,Y\in \cD$ is
  given by
\[
T^\nabla_{X,Y} = \sum_1^3 d\eta_k(X,Y) R_k . 
\]
\end{itemize}
However, since these conditions only place restrictions on the
derivatives in horizontal directions, the connection is not unique and
different extensions are possible.

Now fix a point $p\in Y$. We may assume that the $R_i$, and hence also
the $\eta_i$, are parallel at $p$, i.e.\
\[
\nabla R_i(p)=0, \quad i = 1,2,3.
\]
The connection determines the exponential map in horizontal directions,
\[
\exp_p : \cD_p \longrightarrow Y,
\]
by solving the differential equation $\nabla_{\dot c}{\dot c}=0$,
$c(0)=p$, $\dot{c}(0)=X$, and setting $\exp_p(X)=c(1)$. Linear
coordinates $(y_1,\dots,y_{4m-4})$ on $\cD_p$ give a coordinate system
on the image $S$ of a small ball by $\exp_p$. Then, \emph{at} $p$ one
has
\[
\nabla_{\partial_{y_i}}\partial_{y_j} + \nabla_{\partial_{y_j}}\partial_{y_i} = 
\frac{1}{2} \big(\nabla_{\partial_{y_i}+\partial_{y_j}}(\partial_{y_i}+\partial_{y_j}) - 
\nabla_{\partial_{y_i}-\partial_{y_j}}(\partial_{y_i}-\partial_{y_j}) \big) = 0,
\]
and hence
\[
\nabla_{\partial_{y_i}}\partial_{y_j} = \frac{1}{2} \big( \nabla_{\partial_{y_i}}\partial_{y_j}-
\nabla_{\partial_{y_j}}\partial_{y_i} \big) = \frac{1}{2} T^\nabla_{\partial_{y_i},\partial_{y_j}} 
= \frac{1}{2} \sum_1^3 d\eta_k(\partial_{y_i},\partial_{y_j}) R_k.
\]
In particular, still at $p$,
\[
 \partial_{y_i} \eta_k(\partial_{y_j}) = \eta_k(\nabla_{\partial_{y_i}}\partial_{y_j}) = 
\frac{1}{2} d\eta_k(\partial_{y_i},\partial_{y_j}), 
\]
from which we deduce that for $y\in S$, the projection $Y_i$ of
$\partial_{y_i}$ on $\cD$ satisfies
\[
Y_j = \partial_{y_j} - \frac{1}{2} \sum_1^3 y_i d\eta_k(\partial_{y_i},\partial_{y_j}) R_k + O(|y|^2).
\]
This can be interpreted as saying that if we write the standard
quaternionic contact structure $\Theta_0$ in coordinates $(\sigma_j,y_i)$ on
the Heisenberg group, and denote by $Y_i^0$ the standard horizontal
vector fields as in (\ref{eq:Theta0}) and (\ref{eq:q-left-inv}), then
along $S$ one has
\[
\eta = \Theta_0 + O(|y|^2) , \qquad  Y_i = Y_i^0 + O(|y|^2).
\]
Choosing transverse coordinates $(\sigma_1,\sigma_2,\sigma_3)$ so that
$\partial_{\sigma_i}=R_i$ along $S$, we get the same result in a neighborhood of
$p$ with an error term $O(|y|^2+|\sigma|)$.


\begin{thebibliography}{10}

\bibitem{Bes87}
A.~L. Besse.
\newblock {\em Einstein manifolds}.
\newblock Springer-Verlag, Berlin, 1987.

\bibitem{Biq00}
O.~Biquard.
\newblock Métriques d'{E}instein asymptotiquement symétriques.
\newblock {\em Astérisque}, 265:vi+109, 2000.
\newblock English translation: SMF/AMS Texts and Monographs 13 (2006).

\bibitem{Biq02}
O.~Biquard.
\newblock Métriques autoduales sur la boule.
\newblock {\em Invent. math.}, 148(3):545--607, 2002.

\bibitem{Biq05b}
O.~Biquard, editor.
\newblock {\em Ad{S}/{CFT} correspondence: {E}instein metrics and their
  conformal boundaries}, volume~8 of {\em IRMA Lectures in Mathematics and
  Theoretical Physics}.
\newblock European Mathematical Society (EMS), Zürich, 2005.

\bibitem{ChrDelLeeSki05}
P.~T. Chru{\'s}ciel, E.~Delay, J.~M. Lee, and D.~N. Skinner.
\newblock Boundary regularity of conformally compact {E}instein metrics.
\newblock {\em J. Differential Geom.}, 69(1):111--136, 2005.

\bibitem{Duc04}
D.~Duchemin.
\newblock {\em Géométrie quaternionienne en basses dimensions}.
\newblock PhD thesis, Institut de Recherche Mathématique Avancée,
  Université Louis Pasteur, Strasbourg, 2004.

\bibitem{Duc06}
D.~Duchemin.
\newblock Quaternionic contact structures in dimension 7.
\newblock {\em Ann. Inst. Fourier}, 56(4):851--885, 2006.

\bibitem{EM-tubes}
C.~L. Epstein and R.~B. Melrose.
\newblock {Schrinking tubes and the $\bar{\partial}$-Neumann problem}.
\newblock preprint.

\bibitem{EpsMelMen91}
C.~L. Epstein, R.~B. Melrose, and G.~A. Mendoza.
\newblock Resolvent of the {L}aplacian on strictly pseudoconvex domains.
\newblock {\em Acta Math.}, 167(1-2):1--106, 1991.

\bibitem{Gra00}
C.~R. Graham.
\newblock Volume and area renormalizations for conformally compact {E}instein
  metrics.
\newblock In {\em {The Proceedings of the 19th Winter School ``Geometry and
  Physics'' (Srní, 1999)}}, number~63 in Rend. Circ. Mat. Palermo (2)
  Suppl., pages 31--42, 2000.

\bibitem{GraLee91}
C.~R. Graham and J.~M. Lee.
\newblock Einstein metrics with prescribed conformal infinity on the ball.
\newblock {\em Adv. Math.}, 87(2):186--225, 1991.

\bibitem{Gui05}
C.~Guillarmou.
\newblock Meromorphic properties of the resolvent on asymptotically hyperbolic
  manifolds.
\newblock {\em Duke Math. J.}, 129(1):1--37, 2005.

\bibitem{Helliwell}
D.~Helliwell.
\newblock {Boundary regularity for conformally compact Einstein metrics in even
  dimensions}.
\newblock preprint (based on Thesis (2005) University of Washington), 2006.

\bibitem{Huang}
H.~C. Huang.
\newblock {\em The scattering operator on asymptotically product hyperbolic
  spaces}.
\newblock PhD thesis, Stanford University, 2006.

\bibitem{Lee06}
J.~M. Lee.
\newblock Fredholm operators and {E}instein metrics on conformally compact
  manifolds.
\newblock {\em Mem. Amer. Math. Soc.}, 183(864):vi+83, 2006.

\bibitem{LeeMel82}
J.~M. Lee and R.~Melrose.
\newblock Boundary behaviour of the complex {M}onge-{A}mpère equation.
\newblock {\em Acta Math.}, 148:159--192, 1982.

\bibitem{Mazzeo-EMS}
R.~Mazzeo.
\newblock Analysis and geometry on asymptotically hyperbolic spaces.
\newblock Zürich Lectures in Advanced Mathematics, European Mathematical
  Society.
\newblock Monograph, in preparation.

\bibitem{Maz88}
R.~Mazzeo.
\newblock The {H}odge cohomology of a conformally compact metric.
\newblock {\em J. Differential Geom.}, 28(2):309--339, 1988.

\bibitem{Maz91}
R.~Mazzeo.
\newblock Elliptic theory of differential edge operators. {I}.
\newblock {\em Comm. Partial Differential Equations}, 16(10):1615--1664, 1991.

\bibitem{Maz91b}
R.~Mazzeo.
\newblock Regularity for the singular {Y}amabe problem.
\newblock {\em Indiana Univ. Math. J.}, 40(4):1277--1299, 1991.

\bibitem{MazVas02}
R.~Mazzeo and A.~Vasy.
\newblock Resolvents and {M}artin boundaries of product spaces.
\newblock {\em Geom. Funct. Anal.}, 12(5):1018--1079, 2002.

\bibitem{MazVas04}
R.~Mazzeo and A.~Vasy.
\newblock Analytic continuation of the resolvent of the {L}aplacian on {$\rm
  SL(3)/SO(3)$}.
\newblock {\em Amer. J. Math.}, 126(4):821--844, 2004.

\bibitem{MazVas05}
R.~Mazzeo and A.~Vasy.
\newblock Analytic continuation of the resolvent of the {L}aplacian on
  symmetric spaces of noncompact type.
\newblock {\em J. Funct. Anal.}, 228(2):311--368, 2005.

\bibitem{MazVas07}
R.~Mazzeo and A.~Vasy.
\newblock {Resolvent of the Laplacian on {$\rm SL(3)/SO(3)$}: connections with
  {$3$}-body scattering}.
\newblock {\em Proc. London Math. Soc.}, 2007.

\bibitem{MazMel87}
R.~R. Mazzeo and R.~B. Melrose.
\newblock Meromorphic extension of the resolvent on complete spaces with
  asymptotically constant negative curvature.
\newblock {\em J. Funct. Anal.}, 75(2):260--310, 1987.

\bibitem{Shu01}
M.~A. Shubin.
\newblock {\em Pseudodifferential operators and spectral theory}.
\newblock Springer-Verlag, Berlin, second edition, 2001.
\newblock Translated from the 1978 Russian original by Stig I. Andersson.

\bibitem{Tay96b}
M.~E. Taylor.
\newblock {\em Partial differential equations. {II}}, volume 116 of {\em
  Applied Mathematical Sciences}.
\newblock Springer-Verlag, New York, 1996.
\newblock Qualitative studies of linear equations.

\end{thebibliography}
\end{document}